\documentclass[10pt]{article}
\usepackage{bbm}
\usepackage{eqnarray,amsmath,amsfonts,amsthm,mathrsfs}
\usepackage{color}
\usepackage{bm}
\usepackage{amssymb}
\usepackage{graphicx}
\usepackage{epsfig}
\usepackage{epsf}
\usepackage{float}
\usepackage{subcaption}
\usepackage{amsfonts,amsmath,amsthm,amssymb,graphicx,float,fancyhdr,multirow,hyperref}
\usepackage{booktabs,longtable,authblk}
\usepackage{mathrsfs,hhline}

\usepackage{fancyhdr}
\usepackage[top=2.5cm, bottom=2.5cm, left=3cm, right=3cm]{geometry}
\usepackage{graphicx}

\allowdisplaybreaks[4]

\numberwithin{equation}{section}

\title{Long-time Integration of Nonlinear Wave Equations with Neural Operators$^\dag$\footnotetext{\dag~The work described in this paper is supported partially by Shanghai Science and Technology Program (Project No. 21JC1400600). The work of Zhen Lei is also supported by NSFC Key Program (Grant No. 12431007), the New Cornerstone Science Foundation through the XPLORER PRIZE and Sino-German Center Mobility Programme (Project No. M-0548). The work of Lei Shi is also supported by the NSFC General Program (Grant No. 12171039).}}

\author[1]{Guanhang Lei}
\author[1,2]{Zhen Lei}
\author[1,2]{Lei Shi}

\affil[1]{School of Mathematical Sciences, \linebreak
Shanghai Key Laboratory for Contemporary Applied Mathematics, \linebreak
Fudan University, Shanghai, 200433, China \linebreak
Email: ghlei21@m.fudan.edu.cn \linebreak
\{zlei, leishi\}@fudan.edu.cn}

\affil[2]{Center for Applied Mathematics, \linebreak
Fudan University, Shanghai, 200433, China \linebreak}

\date{}

\begin{document}
	\maketitle
\begin{abstract}
    Neural operators have shown promise in solving many types of Partial Differential Equations (PDEs). They are significantly faster compared to traditional numerical solvers once they have been trained with a certain amount of observed data. However, their numerical performance in solving time-dependent PDEs, particularly in long-time prediction of dynamic systems, still needs improvement. In this paper, we focus on solving the long-time integration of nonlinear wave equations via neural operators by replacing the initial condition with the prediction in a recurrent manner. Given limited observed temporal trajectory data, we utilize some intrinsic features of these nonlinear wave equations, such as conservation laws and well-posedness, to improve the algorithm design and reduce accumulated error. Our numerical experiments examine these improvements in the Korteweg-de Vries (KdV) equation, the sine-Gordon equation, and the Klein-Gordon wave equation on the irregular domain.
\end{abstract}
	

\section{Introduction}\label{section: Introduction}
Solving Partial Differential Equations (PDEs) is central in many physical and engineering fields. Deep learning methodologies and data-driven architectures have become potential alternatives to numerical solvers, considering their low computational costs compared with traditional numerical methods. Among them, neural operators are one class of cutting-edge algorithms. Neural operators learn mappings between infinite-dimensional function spaces, which can be exploited to model the solution operators mapping the parameter functions to solutions, the boundary conditions to solutions, and the initial states to solution trajectories in time-dependent problems. The architectures of neural operators include, but are not limited to, Deep Operator Net (DeepONet) \cite{Lu2021Learning}, Fourier Neural Operator (FNO) \cite{Kovachki2023Neural, Li2021Fourier}, PCA-Net \cite{Bhattacharya2021Model}, Graph Neural Operator (GNO) \cite{Li2020Neural}, multiwavelet-based model (MWT) \cite{Gupta2021Multiwavelet-based}, Wavelet Neural Operator (WNO) \cite{Tripura2023Wavelet}, Operator Transformer (OFormer) \cite{Li2023Transformera}, and General Neural Operator Transformer (GNOT) \cite{Hao2023GNOT}. 

Neural operators have shown advantages compared with conventional methods. Once trained, neural operators can predict the solution given any input function through a network forward calculation process, which is usually a composition of matrix multiplication operations and thus saves computation time compared with traditional numerical methods. Neural operator models are usually meshless methods and resolution-invariant because they can make zero-shot super-resolution predictions: trained on a low-resolution dataset, and give predictions on a high-resolution grid or point cloud. Hence, these models indeed output infinite high-resolution predictions (functions) and learn the operator structure between infinite-dimensional function spaces, which distinguishes them from many traditional methods limited to a mesh structure and constant resolution.

While previous works have applied neural operators to solving various PDE problems, few have considered long-time dynamic system modeling via neural operator methods, especially predicting the long-time integration of nonlinear wave equations. To predict the solution $u$ over a temporal interval $[0, T]$ given initial condition $u_0$, an evident and effective method is training a one-step operator mapping $u_0$ to the complete solution trajectory. However, this method often comes at the price of acquiring long-time observed data. When the length of the data trajectory is significantly insufficient compared to the required prediction time $T$, a natural remedy is to train a short-term predictor and recurrently apply itself to its prediction, which is seen as the next-step initial condition. Although theoretically feasible, many fundamental difficulties will be encountered in numerical implementation. For one-step prediction, the neural network method sacrifices a small amount of accuracy in exchange for a multiple reduction in computation time. This trade-off is attractive and fairly acceptable. However, the error may rise sharply when applying the neural operator to an iterative prediction process. The predictions usually become completely inaccurate, unstable, and blow up in several steps or even the second step. The nonlinearity of the wave equation also introduces nonlinear phenomena, including solitons, shocks, rough waves, peakons, and cnoidal waves. These nonlinear waves usually occur in a long-time evolution, making it difficult to capture by short-term predictors and making this problem more challenging.

However, considering some regularity properties, such as the well-posedness and conservation laws of the nonlinear wave equations, we may expect a numerical stable model that approximately adheres to these properties with linearly growing accumulated error. Regardless of some existing long-time integration models combining neural operators with transfer learning \cite{Xu2023Transfer}, physics-informed neural networks \cite{Wang2023Longtime}, and recurrent neural networks \cite{Michalowska2024Neural}, algorithms or regularization methods inspired by these properties of the nonlinear wave equations themselves are still lacking.

In this paper, we mainly consider FNO, one of the most widely-used neural operators, and examine its numerical performance, particularly its temporal interpolation and extrapolation accuracy, on three types of nonlinear wave equations: Korteweg-de Vries (KdV) equation, sine-Gordon equation, and Klein-Gordon wave equation on an irregular domain. By utilizing some prior knowledge of the equations, such as conservation laws and well-posedness, we propose some approaches that can reduce the accumulation of errors over time. The proposed approaches do not change the general framework designs of neural operators and, hence, can be easily implemented across various neural operator models besides FNO. These include:
\begin{itemize}
\item The use of window-type input data. Considering that the solution of the wave equation is determined by both the initial displacement and velocity, we train neural operators that accept window-type input data composed of multiple displacement snapshots to compensate for the lack of velocity data, eliminating the requirement for collecting velocity observation data.
\item The recurrent neural network architecture. The iterative prediction algorithm aligns with the structure of the recurrent neural network. We embed neural operators into the architecture of the recurrent neural network and compare them with non-recurrent networks. 
\item Random initial conditions distributed on the solution trajectory. We use snapshots at random moments as initial conditions to enhance the representativeness and lower bias of the training data. This approach can significantly improve generalization ability and reduce extrapolated error without increasing training data size. We also discover that it helps models capture the nonlinear phenomena that occur in long-time evolution.
\item Regularization techniques based on the conservation law of energy and momentum. We use the penalty method to impose soft energy constraints on network predictions, serving as a regularization strategy for the optimization.
\item Clipping operator on the prediction. According to the well-posedness and prior Strichartz estimate of the wave equation, we truncate the value of the prediction that exceeds the upper bound. 
\end{itemize}

The remainder of this paper is structured as follows. In \autoref{section: Neural Operators for Nonlinear Wave Equations}, we introduce the long-time integration problem and neural operator models. We also investigate the properties of nonlinear wave equations and suggest some targeted optimization skills for neural operators. In \autoref{section: Numerical Experiments}, we carry out comparative numerical experiments on Korteweg-de Vries (KdV) equation, sine-Gordon equation, and Klein-Gordon wave equation. We demonstrate the numerical experiments process and results. \autoref{section: Conclusions and Discussions} summarizes our findings and discusses future work.

\section{Neural Operators for Nonlinear Wave Equations}\label{section: Neural Operators for Nonlinear Wave Equations}
In this section, we give a brief introduction to neural operators and the long-time integration problem we consider. We draw upon specific insights into the mathematical characteristics of nonlinear wave equations, such as conservation laws and well-posedness, to enhance numerical performance. This allows us to design a more stable and efficient algorithm to utilize the training data and implement neural operators. 

\subsection{Neural operator models}\label{subsection: Neural operator models}
Neural operators are deep learning surrogate models that directly approximate the solution operators of PDEs. These solution operators are mappings between infinite-dimensional function spaces, significantly differing from traditional deep learning models on finite-dimensional tensor spaces. Consider the following family of time-dependent PDEs
\[
    \begin{aligned}
        \mathcal{N}_{a}(u)(t, x) &= f, && (t, x) \in (0, \infty) \times D, \\
        \mathcal{B}(u)(t, x) &= g(t, x), && (t, x) \in (0, \infty) \times \partial D, \\
        \mathcal{I}(u)(t = 0, x) &= u_0(x), && x \in D.
    \end{aligned}
\]
Here, $D \subset \mathbb{R}^d$ is a bounded domain, $a, f, g$ and $u_0$ are functions within specific Banach spaces like $L^p$ or Sobolev spaces. $\mathcal{N}_{a}$ is a parametric differential operator, $\mathcal{B}$ and $\mathcal{I}$ are operators that represent the boundary condition and initial condition, respectively. Solution operators of this PDE can be defined as $\mathcal{G}^{\dagger}: (a, f, g, u_0) \mapsto u$, or any component of it such as $\mathcal{G}^{\dagger}: a \mapsto u$ when other input functions are fixed. The domain of definition and range of $\mathcal{G}^{\dagger}$ are the corresponding Banach spaces, say $\mathcal{G}^{\dagger}: \mathcal{A} \to \mathcal{U}$. Neural operators build an approximation of $\mathcal{G}^{\dagger}$ by constructing a parametric map $\mathcal{G}_{\theta}: \mathcal{A} \to \mathcal{U}$ with $\theta \in \mathbb{R}^{p}$. This process is usually done by training a parametric neural network with some observation data $\{a_i, u_i\}_{i=1}^N$ through a supervised learning paradigm, minimizing the empirical risk functions measuring the error between the prediction $\mathcal{G}_{\theta}(a_i)$ and ground truth $u_i$. Once the training process is finished, the neural operator $\mathcal{G}_{\theta}$ is ready to give prediction $\mathcal{G}_{\theta}(a)$ given any input function $a$. Hence, it is an efficient model using an offline-online computational strategy that not only solves one specific PDE, but a family of them taking any possible input function. We also note that some unsupervised operator learning methods do not require any paired input-output data $(a, u)$ by assuming that the explicit form of the equation is known. These are frequently referred to as the physics-informed models, such as physics-informed neural network (PINN) \cite{Sirignano2018DGM, Raissi2019Physicsinformed}, physics-informed DeepONet \cite{Wang2021Learning} and physics-informed FNO \cite{Li2024Physics-Informed}. They train the models by minimizing the associated PDE residuals $\vert \mathcal{N}_{a}(u) - f \vert, \vert \mathcal{B}(u) - g \vert$ and $\vert \mathcal{I}(u) - u_0 \vert$. Our discussion is restricted to a data-driven regime, so we will not consider physics-informed models in the following sections. 

DeepONet \cite{Lu2021Learning} and FNO \cite{Kovachki2023Neural, Li2021Fourier} are two popular neural operator architectures. DeepONet is inspired by the operator networks proposed in \cite{TianpingChen1995Universal}, where they proved that these networks possess a universal approximation property for infinite-dimensional nonlinear operators. A DeepONet consists of two deep neural networks, termed branch net and trunk net. The branch net takes the input function $a$, encodes it as a finite-dimensional vector by collecting its evaluated values at a set of fixed sensors $\{x_i\}_{i=1}^m \subset D$, and extracts latent representations through a deep feedforward neural network. The trunk net takes the coordinates $y$ where the output functions are evaluated and also extracts latent representations through a neural network. Finally, the query value $\mathcal{G}_{\theta}(a)(y)$ is computed by a dot product of the outputs of branch net and trunk net, similar to the process of multiplying basis functions by coefficients. FNO was initially proposed to parameterize the integral kernel operator, which is a more general neural operator model based on the form of the Green's function solution. FNO lifts the input function to a higher dimensional latent space, applies layers of integral kernel operators, and then projects back to the solution space. In each layer of integral kernel operators, the latent representations go through a fast Fourier transform (FFT), a linear transform on only the lower Fourier modes, an inverse FFT, a weighted residual link, and finally the nonlinear activation. 

Both DeepONet and FNO are resolution-invariant in that they can be trained and evaluated on different resolution data and even on various domains of definition. DeepONet is a meshless method that only takes the coordinates of query points as input for the trunk net. However, the branch net only accepts input of a fixed resolution as it directly applies a fully-connected neural network on the point-wise evaluated values. When encountering input function data with a different resolution, the DeepONet needs to be retrained. The BasisONet model proposed in \cite{Hua2023Basis} uses numerical integrations as the resolution-invariant function autoencoder that can work on different discretization inputs. FNO can be trained and evaluated on different resolution data but is restricted initially to rectangular grids required by FFT. Some subsequent works successfully extended the FNO to irregular geometries. Geo-FNO in \cite{Li2023Fourier} uses a trainable deformation to transform the irregular mesh on a general physical geometry to a rectangular grid on a latent computational space. \cite{Liu2023Domain} encodes the geometry information by incorporating the domain characteristic function into the integral layer of FNO. \cite{Li2023Geometry-Informed} combines GNO and FNO architectures, thereby leveraging the ability of GNO to handle irregular grids via graph-based methods. In \autoref{subsection: Klein-Gordon wave equation on irregular domain}, we will consider the Klein-Gordon wave equation on an irregular domain using Geo-FNO.

\subsection{The long-time integration problem}\label{subsection: The long-time integration problem}
In this paper, we aim to predict the long-time integration of nonlinear wave equations. Specifically, we consider the following equation
\begin{align}
    \partial^2_t u(t, x) - \Delta u(t, x) &= f, && (t, x) \in (0, \infty) \times D, \label{differential equation}\\
    \mathcal{B}(u)(t, x) &= g(t, x), && (t, x) \in (0, \infty) \times \partial D, \label{boundary condition}\\
    u(t = 0, x) &= u_0(x), && x \in D, \label{initial displacement}\\
    \partial_t u(t = 0, x) &= v_0(x) && x \in D. \label{initial velocity}
\end{align}
This equation is of rich physical implications, where $u$ represents the displacement scalar field and $\partial_t u$ represents the wave velocity. As we consider nonlinear cases, $f$ can be dependent on $u$ and its high-order derivatives. Assume that the equation \eqref{differential equation} and boundary condition \eqref{boundary condition} are fixed, we learn the solution trajectory over a long-time temporal interval $[0, T]$ from given initial conditions input pairs $(u_0, v_0)$. If there is sufficient available complete trajectory data, one could immediately train the neural operator approximating $(u_0, v_0) \mapsto u\vert_{t \in [0, T]}$. The essential difficulty of this problem is the "long-time" mentioned here, which means that the data length we can use to train is very limited compared with the required prediction time $T$. 

The most natural idea to address this issue, which is also adopted in \cite{Tripura2023Wavelet, Wang2023Longtime, Michalowska2024Neural, Xu2023Transfer, Navaneeth2024Waveformer}, is to train a short-time propagator $\mathcal{G}: (u_0, v_0) \mapsto u\vert_{t \in [0, \Delta T]}$ using short-term data, and in each step, recurrently applied $\mathcal{G}$ to its predicted solution at the final stage, until the predicted time $[0, n \Delta T]$ fully covers $[0, T]$ after $n$ steps of recursion. Here $\Delta T$ is usually much smaller than $T$, generally resulting in a large $n$. As one might expect, the predictions and accumulated error usually blow up as $n$ rises. In many cases, the predictions become imprecise even in the second step, indicating that the model has bad extrapolated accuracy and generalization ability outside the training region. This failure is mainly attributed to the insufficient prediction of nonlinear phenomena, including solitons, shocks, and rough waves. A soliton does not change shape when it propagates, even after interacting with other solitons. Rough waves and shocks are huge, steep waves that emerge unpredictably from smooth waves as time evolves. Some nonlinear phenomena only occur when the wave has enough time to evolve. Hence, their mechanisms are complex for short-term predictors to grasp. If only the solutions from a limited period are used as the training set, the network may fail to predict the generation of nonlinear waves. This also highlights the difficulties in long-time integration problems. Now, we discuss several methods that may enhance the generalization ability of neural operators.

\textbf{Temporal window input:} 
Some recent works also employ this recurrent idea in other time-dependent equations with a first-order partial derivative with respect to time, such as the KdV equation \cite{Wang2023Longtime, Michalowska2024Neural}, Allen-Cahn equation \cite{Xu2023Transfer, Navaneeth2024Waveformer} and reaction-diffusion equation \cite{Wang2023Longtime, Xu2023Transfer}. All of them only consider using one predicted snapshot as input to the next step prediction, that is, for $1 \leq i \leq n-1$,
\[
u \vert_{t \in (i\Delta T, (i+1)\Delta T]} = \mathcal{G}_{\theta}(u\vert_{t = i\Delta T}).
\] 
Although initial displacement is enough to determine the solution of these equations, it is not feasible in the wave equation with second-order time partial derivative as the initial velocity also affects the solution, suggesting that one should use the following recursion:
\[
\begin{pmatrix}
u \vert_{t \in (i\Delta T, (i+1)\Delta T]} \\
\partial_t u \vert_{t \in (i\Delta T, (i+1)\Delta T]}
\end{pmatrix}
=
\mathcal{G}_{\theta}
\begin{pmatrix}
u \vert_{t = i\Delta T} \\
\partial_t u\vert_{t = i\Delta T}
\end{pmatrix}.
\]
This puts higher demands on both data and the model. One has to collect the velocity trajectory data and train a larger-scale network to deal with additional velocity information. In reality, it can be expensive to track the velocity trajectory. The neural operator also needs to figure out the coupling dynamic between displacement and velocity data with inconsistent distributions, making the training process more costly. 

Alternatively, to collect the velocity data, we can use the forward-difference method or high-order methods on the displacement data to numerically approximate the velocity. Then, we train the network with concatenated displacement-velocity data. As an alternative, we can also use displacement data over a temporal window as input: 
\begin{equation}\label{window input}
\mathcal{G}_{\theta}: u\vert_{t \in [(i-1)\Delta T, i\Delta T)} \mapsto u\vert_{t \in [i \Delta T, (i+1) \Delta T)}.
\end{equation}
Data in this temporal window is of the time-discretized form, for example, the uniform discretization 
\begin{equation}\label{window discretization}
u\vert_{t \in [(i-1)\Delta T, i\Delta T)} = u\vert_{t \in \{(i-1)\Delta T, (i-1)\Delta T + \Delta t, \ldots, (i-1)\Delta T + (l-1) \Delta t\}},
\end{equation}
where $l = \Delta T/\Delta t$ represents the width of the temporal window. We may expect that the network can automatically extract velocity features or any other useful features (not necessarily the velocity) from this temporal input window that help predict the trajectory of the solution. Roughly speaking, for the wave equation, the initial displacement in the fine-discretized window also uniquely determines the solution under the condition of modulating some high-frequency waves. Using window form input has several advantages compared with paired displacement-velocity data, which eliminates the necessity to select a numerical differentiation method and avoids the error caused by approximating velocity. Besides, the distribution of data shares a more similar pattern, and the network is an end-to-end architecture that is more suitable for both training and generalization.

\begin{figure}[t]
    \centering
    \includegraphics[width=\textwidth]{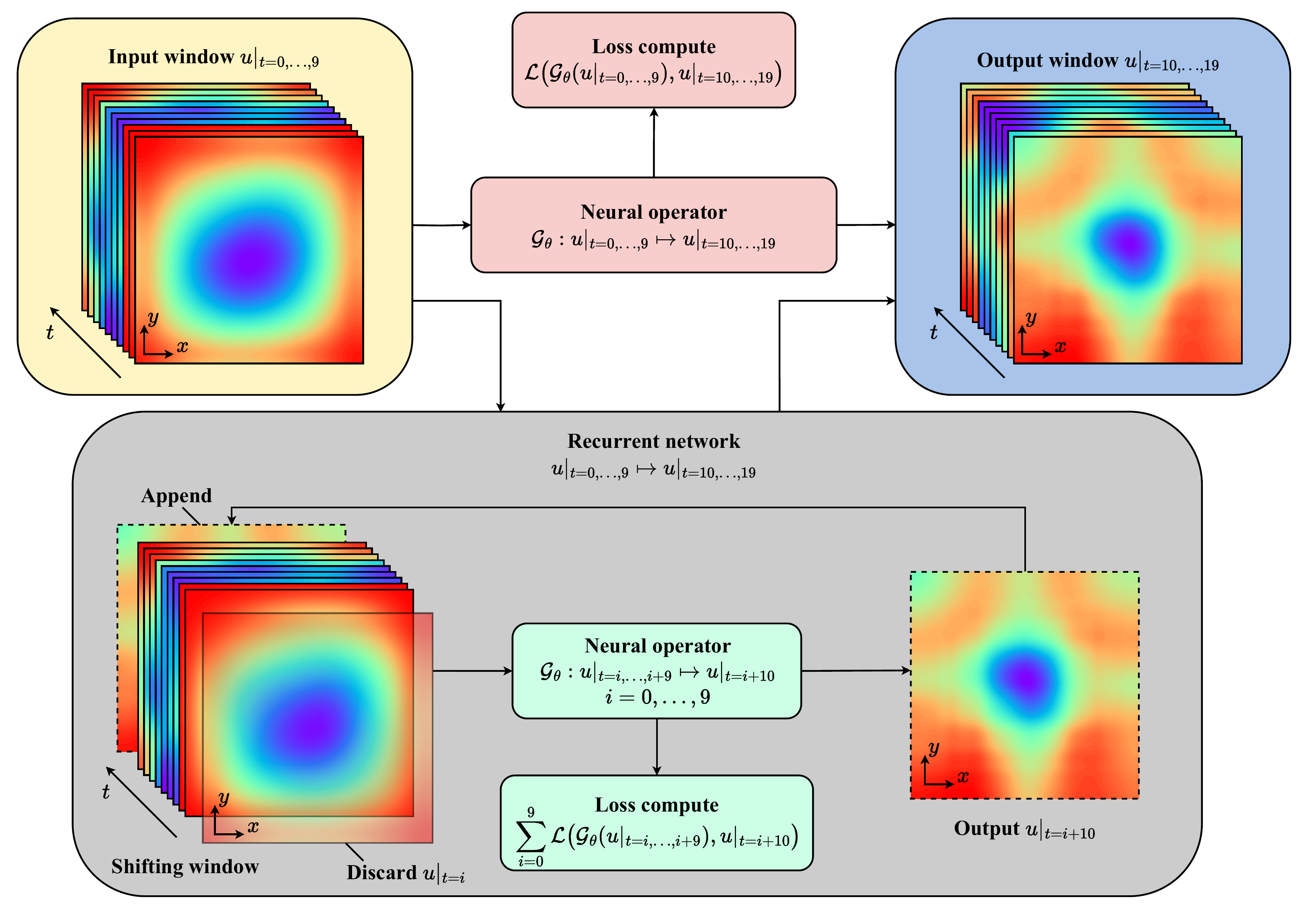}
    \caption{\textbf{The recurrent network and full prediction network.} The operator to be learned exhibits a window-to-window mapping structure. The full prediction network (red box) learns this mapping directly and predicts the full output window in one step. The recurrent network (grey box) uses the neural operator to predict a snapshot of the next moment in each step. The snapshot is then appended to the shifting input window and the first moment snapshot is discarded. This process continues until the next window is completely predicted.}
    \label{figure:window}
\end{figure}

\textbf{Recurrent network and full prediction network:} 
We propose two architectures of the network implementing the recursion \eqref{window input}, see \autoref{figure:window}. The first one adopts a recurrent neural network (RNN) structure within the operator. The RNN takes the input window \eqref{window discretization} and predicts the solution at $t = i\Delta T$. This new prediction is then appended to the back end of the window, and the input data of the first moment $t = (i-1)\Delta T$ is discarded. This process continues until the next window prediction $u\vert_{t \in \{i\Delta T, i\Delta T + \Delta t, \ldots, i\Delta T + (l-1) \Delta t\}}$ is completed and the original input is fully discarded. The process is similar to shifting the window on the timeline step by step. The next window prediction will be compared with the training data to calculate the loss. Hence, each prediction step contributes to the loss, and the gradient backpropagation will occur in the stacked recurrent network structure. The second architecture directly learns the window-to-window mapping \eqref{window input}, making the full prediction in each step. In \eqref{window input}, input and output share the same temporal length. One may also consider using different window lengths between input and output. 

The recurrent operator model requires more iteration steps during inference, while its error in each step may be lower as its output shape is simpler. The recurrent training method may be more suitable for recurrent prediction as they share the same iterative manner. We can also leverage various well-developed techniques and structures for recurrent neural networks to improve the recurrent operator model, such as long short-term memory (LSTM), gated recurrent units (GRU), gradient clipping, attention mechanism, etc. We will compare the performances of recurrent and full prediction architectures in numerical experiments. 

\textbf{Random initial time:}
In our numerical experiments, we generate the data set by solving the nonlinear wave equation using traditional numerical solvers with $v_0 = 0$ and $u_0$ sampled from a distribution, such as the Gaussian random field. We obtain the entire solution trajectory in a long temporal interval $[0, T]$, and use truncated parts $u \vert_{t \in [0, \Delta T)}$ as input and $u \vert_{t \in [\Delta T, 2\Delta T)}$ as output to train the neural operator model \eqref{window input}, where $i = 1$. We take $\Delta T$ much smaller than $T$, for instance, $\Delta T = T/10$. Hence, the size of the training set is very limited, while we require the model to predict the complete solution $u \vert_{t \in [0, T]}$. Although we generate the initial displacement in the Gaussian random field, which would be regarded as a "sufficiently random" function, the distributions of the training set and the test set are inconsistent in this task: When training the neural operator, the model consistently receives input functions from the starting part of a trajectory, where the velocity field is approximately zero. However, during iterative prediction, the model has to deal with input functions that are distributed across the entire trajectory. This discrepancy introduces data bias and causes low generalization ability in long-time prediction.

To address this issue, we choose a random initial time $t_i \in [0, T - 2 \Delta T]$ for each trajectory sample $u_i \vert_{t \in [0, T]}$ and use the truncated part $(u_i \vert_{t \in [t_i, t_i + \Delta T)}, u_i \vert_{t \in [t_i + \Delta T, t_i + 2\Delta T)})$ to train the network. Even though this does not increase the size of the training set and the operator we train remains a short-term prediction model, we find it very effective in reducing accumulated error, as it lowers the bias of training data distribution among the entire trajectory. This data randomness is often more important than the randomness introduced by generating initial conditions from a random function field in reducing accumulation error, as it allows the model to learn nonlinear phenomena, such as shocks and rough waves, that only appear after long-term evolution.

We choose only one fragment from every trajectory for training. It is also possible that all trajectories are divided into fragments, and all fragments are used for training. There are two primary reasons for our paper not to consider this method. First, if all fragments from a trajectory are available, one can concatenate them to form the complete trajectory and train a long-term predictor directly, which eliminates the necessity of iterative prediction and usually enjoys better stability. Second, using all fragments of all trajectories significantly increases the data size and training time.

\subsection{Conservation laws and well-posedness of nonlinear wave equations}\label{subsection: Conservation laws and well-posedness of nonlinear wave equations}
We now discuss some intrinsic properties of nonlinear wave equations, which inspire our algorithm designs for neural operators. These topics have central roles in the PDE theory and have been extensively studied in various cases, including semilinear equations, quasilinear equations, the Cauchy problem, different boundary conditions, and cases involving different spatial dimensions and regularity function spaces. As it is impossible to cover general cases, we will pinpoint specific examples in our discussion.

\textbf{Conservation laws:}
We first consider the following KdV equation 
\[
\partial_t u + u \partial_x u + \partial^3_x u = 0, \quad (t,x) \in (0, T] \times \mathbb{R}.
\]
Though it is not of the form \eqref{differential equation}, it usually models shallow water waves, acoustic waves, and long internal ocean waves \cite{Miles1981Korteweg-Vries} and shares many properties with the wave equation \eqref{differential equation}. Hence we adopt this equation as our first numerical example, which also aligns with several previous works \cite{Wang2023Longtime, Michalowska2024Neural} on neural operators. It is well-studied \cite{Miura1968Korteweg-Vries} that the KdV equation has an infinite number of conserved quantities that do not change in time and we list the first few explicitly here:
\[
E_1(u) = \int^{+\infty}_{-\infty} u \mathrm{d}x, \quad E_2(u) = \int^{+\infty}_{-\infty} u^2 \mathrm{d}x, \quad E_3(u) = \int^{+\infty}_{-\infty} \frac{1}{3}u^3 - (\partial_x u)^2 \mathrm{d}x.
\]
All these conservation laws are integrals of polynomials with respect to $u$ and its various order spatial derivatives. These conservation laws are given in the case of the one-dimensional Cauchy problem. One can easily generalize them to initial-boundary value problems including the one-dimensional periodic condition problem in our numerical example. 

Considering that predictions by neural operators usually become unstable and blow up in several iteration steps, we impose soft constraints on the output of the operator, by adding the difference between the conserved quantities of input and output to the loss function:
\begin{equation}\label{conservation law loss}
\mathcal{L}(\theta) = \mathcal{L}_0\big(\mathcal{G}_\theta(u_{\mathrm{in}}), u_{\mathrm{out}}\big) + \lambda_1 \mathcal{L}_{\mathrm{cl}}\big(E_1(\mathcal{G}_\theta(u_{\mathrm{in}})), E_1(u_{\mathrm{out}})\big) + \lambda_2 \mathcal{L}_{\mathrm{cl}}\big(E_2(\mathcal{G}_\theta(u_{\mathrm{in}})), E_2(u_{\mathrm{out}})\big) + \cdots,
\end{equation}
where $\mathcal{L}_0$ is the loss function used to fit the data, $\mathcal{L}_{\mathrm{cl}}$ measures the model's deviation on the conservation law, and $\lambda_i$ are regularization hyperparameters. For example, we can take the square loss $\mathcal{L}_{\mathrm{cl}}(x,y) = (x-y)^2$. In this case, we only use the first two conservation laws $E_1, E_2$ as they are easily calculated, while the third one involves the calculation of the derivative. The conserved quantities serve as regularization techniques and we expect the reduction of the large-scale deviations and instability of the output.

We then consider the following nonlinear wave equation
\[
\partial^2_t u(t, x) - \Delta u(t, x) = F'(u), \quad (t, x) \in (0, T] \times D
\]
subject to homogeneous Dirichlet or Neumann boundary conditions with a potential force function $F(u)$. It is also well-known that this equation is energy-conserving:
\[
\frac{\mathrm{d}}{\mathrm{d}t}E(t) := \frac{\mathrm{d}}{\mathrm{d}t} \int_D (\partial_t u(t, x))^2 + \|\nabla u(t, x)\|^2 + 2F(u(t, x)) \mathrm{d}x = 0.
\]
Our conservation quantities regularization technique can also be applied in this case, by some additional calculations on the derivatives. One simple way to approximate the derivatives is to use the difference method on the discretized output, which may fail when the discretization is on irregular mesh or even meshless point cloud. An alternative way is to use the autograd toolkit given that our model is built upon feedforward networks, which is similar to the derivative calculation process in physics-informed neural networks. We also note that this technique is based on prior energy-conserving training data, which requires energy-conserving numerical solvers (see \cite{Brugnano2015Energy, Li2020Linearly}).

Our soft constraint method originates from the penalty method in optimization. For the quadratic penalty method, it is demonstrated that when the penalty parameter $\lambda \to \infty$, the solution of the optimization problem will converge to a KKT point or global minimum point, which indicates that a perfect-optimized model of our problem is energy-preserving. This inspires us to use a dynamically increasing regularization parameter, multiplying $\lambda$ by some $a > 1$ at each optimization step.

\textbf{Well-posedness:} 
A PDE is called well-posed in the sense of Hadamard if a solution exists, the solution is unique and depends continuously on the conditions. We focus on the continuity here, especially the continuity with respect to initial values. If continuous dependence on the data can be established, we can theoretically expect the stability of the operators under minor perturbations and good generalization on a similar dataset. Energy estimates are critical tools to establish the well-posedness of wave equations. Abundant results have been established for both Cauchy problems \cite{Tao2006Nonlinear} and initial-boundary problems. The so-called Strichartz estimates are typical energy estimates for dispersive PDEs. We quote the following Strichartz estimates for the nonlinear wave equation \eqref{differential equation}\eqref{boundary condition}\eqref{initial displacement}\eqref{initial velocity}. Readers may refer to Corollary 1.2. of \cite{Blair2009Strichartz}, Proposition 3.1. of \cite{Burq2008Global}, and Proposition 2.1. of \cite{Burq2009Global} for details regarding the conditions and parameters.
\[
\begin{aligned}
&\|u\|_{L^p([-T, T]; L^q(M))} \leq C\big( \|u_0\|_{H^\gamma(M)} + \|v_0\|_{H^{\gamma-1}(M)} + \|f\|_{L^r([-T, T]; L^s(M))} \big), \\
&\|u\|_{L^5((0, 1); W^{\frac{3}{10}, 5}_0(D))} + \|u\|_{C^0((0, 1); H^1_0(D))} + \|\partial_t u\|_{C^0((0, 1); L^2(D))} \\
&\quad \quad \leq C\big( \|u_0\|_{H^1_0(D)} + \|v_0\|_{L^2(D)} + \|f\|_{L^\frac{5}{4}((0, 1); W^{\frac{7}{10}, \frac{5}{4}}(D))} \big).
\end{aligned}
\]
Strichartz estimates give us a priori bound on the solution, which can be seen as a supplement of the conservation law and hence another regularization technique for network training. We can also apply a clipping operation on prediction by truncating the value that exceeds the upper bound:
\[
\pi(x):=
\begin{cases}
-B & x<-B \\
t & x \in [-B,B] \\
B & x>B.
\end{cases}
\]
The upper bound $B$ can be chosen by calculating the Strichartz estimates above with a suitable $C$ or selecting an empirical suitable $B$ according to the training dataset. The truncating function $\pi(x)$ can be only used in the predicting process, since during the training process, the truncation will overwrite the loss information that should participate in gradient backpropagation. To use the clipping operator in the training process, one should use a soft truncation function, like $\sigma(x) = B\cdot\tanh(x)$, as the activation function added right before the output of the neural operator. 

Clipping operation is widely used in various machine learning algorithms, including neural networks \cite{Bartlett1998sample} and support vector machines \cite{Bousquet2002Stability} to enhance stability and generalization. Clipping operation also guarantees the boundedness of the prediction, thereby facilitating theoretical analysis and the derivation of oracle inequalities. We emphasize that the well-posedness and these estimates are crucial conditions to establish the fast convergence rate on the generalization upper bound of the PDE learning algorithm. In our previous work \cite{Lei2025Solving}, we derive rigorous assumptions on the well-posedness and regularity of PDEs to establish the strong convexity of the PINN risk with respect to the Sobolev norm. This convexity is crucial to bound the approximation and generalization error and lead to a fast learning rate of the algorithm, which is also mentioned in related works \cite{Lu2022Machine, Jiao2022rate, Mishra2023Estimates}.

\section{Numerical Experiments}\label{section: Numerical Experiments}
In this section, we conduct numerical experiments to verify our discussion in \autoref{section: Neural Operators for Nonlinear Wave Equations}. We consider the ($1+1$)-dimensional KdV equation with periodic boundary condition, ($2+1$)-dimensional sine-Gordon equation with Neumann boundary condition, and ($2+1$)-dimensional Klein-Gordon wave equation on an irregular domain with Dirichlet boundary condition. Because some of the properties we discussed before may apply to specific types of equations, we only compare the numerical performance for certain methods in each example. Most basic settings between different experiments are the same, which we now explain here as follows.

In all cases, we generate the data using traditional PDE solvers. The complete trajectory of each data is uniformly discretized to at least 100 snapshots. We only use 20 continuous snapshots as the training set (always the first 20 snapshots when not using random initial time). We use temporal windows to train the neural operator to learn \eqref{window input}. Input and output windows both cover 10 snapshots, hence the window size $l = 10$. A comparison on different window sizes $l$ is considered in \autoref{subsection: Sine-Gordon equation}. The structure of the model can be either a recurrent network or a full prediction network. For optimization, we always use ADAM optimizer with an initial learning rate of 0.001 with weight decay of $10^{-6}$. The learning rate decays at a rate of 0.75 every 50 epochs. We split a validation set from the training set and choose the model with the smallest validation error. The data size, batch size, and epochs are listed in \autoref{tab:Data size, batch size and number of epochs}. At the inference stage, given the initial 10 snapshots, the model predicts the complete trajectory of the testing data. We calculate the average $L^2$ relative error of the predictions across the test set at each moment and plot the accumulated error. All numerical experiments are performed on a single Tesla P100-PCIE-16GB GPU. 

\begin{table}[t]
    \centering
    \caption{\textbf{Data size, batch size and number of epochs.}}
    \begin{tabular}{lccccc} 
        \toprule
        PDE & Training size & Validation size & Test size & Batch size & Number of epochs \\ 
        \midrule
        KdV & 800 & 200 & 200 & 5 & 500 \\
        Sine-Gordon & 400 & 100 & 100 & 10 & 1000 \\
        Klein-Gordon & 400 & 100 & 100 & 10 & 1000 \\
        \bottomrule
    \end{tabular}
    \label{tab:Data size, batch size and number of epochs}
\end{table}
\begin{table}[t]
    \centering
    \caption{\textbf{Number of parameters and time per epoch.}}
    \begin{tabular}{lcccccc}
        \toprule
        \multirow{2}{*}{Model} & \multicolumn{3}{c}{Number of parameters} & \multicolumn{3}{c}{Time per epoch(s)} \\
        \cline{2-7}
                                & KdV & SG & KG & KdV & SG & KG \\
        \midrule
        Recurrent DeepONet & 223,409 & - & - & 8.91 & - & - \\
        Full prediction DeepONet & 223,509 & - & - & 4.15 & - & - \\
        Recurrent FNO & 255,745 & 2,061,739 & 23,615,681 & 30.19 & 16.29 & 45.37 \\
        Full prediction FNO & 299,393 & 2,101,649 & 23,616,842 & 7.26 & 7.94 & 5.03 \\
        \bottomrule
    \end{tabular}
    \label{tab:Number of parameters and time per epoch}
\end{table}

\subsection{KdV equation}\label{subsection: KdV equation}
As we have mentioned in \autoref{section: Neural Operators for Nonlinear Wave Equations}, regardless of its different form from \eqref{differential equation}, the KdV equation was initially used to model shallow water waves and also describes the evolution of various physical waves such as nonlinear acoustic wave, gravity waves, and plasma \cite{Miles1981Korteweg-Vries}. The solution of the KdV equation usually exhibits solitons, shocks, and rough wave phenomena, which introduces more significant challenges to prediction.

Here, we consider the following KdV equation
\[
\begin{aligned}
\partial_t u + u \partial_x u + \delta^2 \partial^3_x u &= 0, \quad \delta = 0.01, && (t,x) \in (0,1) \times (0,1), \\
u(t, x = 0) &= u(t, x = 1), && t \in [0, 1], \\
u(t = 0, x) &= u_0(x), && x \in [0, 1].
\end{aligned}
\]
We generate 1200 trajectories with 1000 training samples and 200 testing samples. Each trajectory is uniformly discretized with a spatial resolution of 1024 and 100 temporal snapshots including the first snapshots $u_0$. The initial condition $u_0$ is drawn from the Gaussian random field $\mathcal{N}(0, 7^4(-\Delta + 7^2I)^{-2.5})$ with periodic boundary conditions, and the equation is solved using chebfun package \cite{Platte2010Chebfun}, aligned with \cite{Gupta2021Multiwavelet-based}. 

As the first trial, we test the performance of DeepONet and FNO, both using recurrent network and full prediction network structures. The results presented in \autoref{tab:KdV interpolation error and extrapolation error} and \autoref{fig:KDV-draw_error_DON-FNO} demonstrate that FNO significantly outperforms DeepONet. Hence, in subsequent experiments, we only consider the FNO architecture. We note that the FNO in a recurrent network uses 1-dimensional convolution in each latent layer, i.e. the FNO-1D model, while the FNO in a full prediction network is FNO-2D. This is because the input (window) and output (snapshot) shapes of the recurrent network are different. Hence, we need to view the temporal dimension as the channel dimension, and only the spatial dimension goes through the convolution operator. A full prediction network has the same input (window) and output (window) shape so 2D convolution is allowed. 

We have discussed the infinite conservation laws and designed the regularization loss \eqref{conservation law loss} for the KdV equation in \autoref{section: Neural Operators for Nonlinear Wave Equations}. \autoref{tab:KdV interpolation error and extrapolation error}, \autoref{fig:KDV-draw_error_Win-ConLaw} and \autoref{fig:KDV-draw_error_Block-ConLaw} compares FNO trained with regularized loss \eqref{conservation law loss} setting $\lambda_1 = \lambda_2 = \lambda$ and vanilla FNO ($\lambda = 0$). We conclude that a suitable chosen $\lambda$ can significantly reduce the extrapolation error. We illustrate the prediction generated by the full prediction FNO with $\lambda = 10$ in \autoref{fig:KDV-plot-prediction}, which has a relatively low absolute error mainly concentrated around the rough wave (the red stripe at $t > 0.6$ in \autoref{fig:KDV-plot-prediction}). By contrast, the prediction of soliton (the blue stripe) is more accurate.

To address the issue of nonlinear phenomena prediction, we then use random initial time data to train the model, which is expected to better simulate the generation of rough waves as time evolves, compared with models trained on vanilla data. We consider two ways to generate the random initial time. The first is to directly generate moments that are uniformly distributed on the entire admissible interval $t \in [0, 0.80]$, we call it global random data. The second focuses on the time around the occurrence of rough waves. For this KdV problem, we consider the latter half of the time period: $t \in [0.50, 0.80]$ and we call it local random data. We show the accumulation error of models trained with vanilla data, global random data, and local random data in \autoref{tab:KdV interpolation error and extrapolation error}, \autoref{fig:KDV-draw_error_Win-RandomStart} and \autoref{fig:KDV-draw_error_Block-RandomStart}. Numerical results demonstrate that using random data, especially local random data, can trade a degree of short-term prediction accuracy for significantly improved stability in long-term prediction. We draw the prediction made by the full prediction FNO trained with local random data in \autoref{fig:KDV-plot-prediction2} and discover that its absolute error is low and evenly distributed in space, which indicates that the model can also predict spatial-localized nonlinear waves well.

\begin{table}[t]
    \centering
    \caption{\textbf{Interpolation error and extrapolation error for KdV equation.}}
    \begin{tabular}{lcc} 
        \toprule
        \multirow{2}{*}{Models} & Interpolation & Extrapolation \\ 
        & $t \in [0.10, 0.19]$ & $t \in [0.20, 0.99]$ \\
        \midrule
        Recurrent DeepONet  & 0.1266 & 0.7166 \\
        Full prediction DeepONet & 0.0988 & 0.7198 \\
        \midrule
        Recurrent FNO, $\lambda = 0$ & 0.0081 & 0.3416 \\
        Recurrent FNO, $\lambda = 0.1$ & 0.0080 & \textbf{0.2198} \\
        Recurrent FNO, $\lambda = 1$ & 0.0083 & 0.4049 \\
        Recurrent FNO, $\lambda = 10$ & 0.0081 & 0.3892 \\
        Recurrent FNO, $\lambda = 100$ & \textbf{0.0078} & 0.3639 \\
        \midrule
        Full prediction FNO, $\lambda = 0$ & 0.0128 & 0.2544 \\
        Full prediction FNO, $\lambda = 0.1$ & \textbf{0.0122} & 0.2492 \\
        Full prediction FNO, $\lambda = 1$ & 0.0125 & 0.2381 \\
        Full prediction FNO, $\lambda = 10$ & 0.0126 & \textbf{0.2275} \\
        Full prediction FNO, $\lambda = 100$ & 0.0123 & 0.2525 \\
        \midrule 
        Recurrent FNO, vanilla data & \textbf{0.0081} & 0.3416 \\
        Recurrent FNO, global random data & 0.0095 & 0.1820 \\
        Recurrent FNO, local random data & 0.0139 & \textbf{0.1023} \\
        \midrule 
        Full prediction FNO, vanilla data & 0.0128 & 0.2544 \\
        Full prediction FNO, global random data & \textbf{0.0126} & 0.2080 \\
        Full prediction FNO, local random data & 0.0227 & \textbf{0.0802} \\
        \bottomrule
    \end{tabular}
    \label{tab:KdV interpolation error and extrapolation error}
\end{table}

\subsection{Sine-Gordon equation}\label{subsection: Sine-Gordon equation}
The sine-Gordon equation is a nonlinear dynamical model introduced initially in differential geometry. It is Lorentz invariance and arises in many physical models, including nonlinear optics, DNA-soliton dynamics, crystal dislocation, nonlinear theories of elementary particles, etc. Readers can refer to \cite{Barone1971Theory, Cuevas-Maraver2014sine-Gordon} for details.

We consider a ($2+1$)-dimensional case with a homogeneous Neumann boundary:
\[
    \begin{aligned}
        \partial_t^2 u - \Delta u + \sin u &= 0, && (t, x) \in (0, 20) \times (0, 1)^2 ,\\
        \frac{\partial u}{\partial n} (t, x) &= 0 , && t \in (0, 20), x \in \partial (0, 1)^2, \\
        u(t = 0, x) &= u_0(x), && x \in (0, 1)^2, \\
        \partial_t u(t = 0, x) &= 0, && x \in (0, 1)^2. 
    \end{aligned}
\]
The data set contains $600$ samples with $100$ testing samples. The spatial-temporal resolution is $64^2 \times 200$ with the first snapshot $u_0(x)$ drawn from the Gaussian random field $\mathcal{N}(0, 10^4(-\Delta + 8^2I)^{-6})$ with Neumann condition. The PDE is solved using \cite{Su2019Numerical}, a localized method of approximate particular solutions with radial basis function finite difference method. 

We use global random initial time data to train the networks and compare them with the vanilla data using original initial conditions. \autoref{tab:SG interpolation error and extrapolation error} and \autoref{fig:SG-draw_error} illustrate the performance of recurrent FNO (FNO-2D), full prediction FNO (FNO-3D) using vanilla data and random data. For models trained with vanilla data, the full prediction FNO blows up immediately in the extrapolation step, and the recurrent FNO has an error that oscillates upwards periodically, due to the periodicity of the wave. When the wave returns to a phase similar to the initial value, the operator shows a better generalization performance and the prediction error will significantly decrease. The numerical results indicate that, despite the reasonable sacrifice of interpolation accuracy, using random initial time data can suppress long-time error fluctuation and improve generalization as the model can learn the behavior of the wave at different phases. This can be seen as a trade-off between the bias and variance of the model. \autoref{fig:SG-plot-prediction} shows the prediction given by recurrent FNO trained with random initial data. 

This paper considers window-type input data as a substitute for displacement-velocity data. The window size $l$ also significantly impacts model performance. To show this point, we train the recurrent FNO with random data and different window sizes $l$ for the sine-Gordon equation and plot their accumulation error in \autoref{fig:SG-windowsize}. Their average relative errors are shown in \autoref{tab:SG-l}. We find that the model with $l=1$ (not using window-type input data) does not converge, which corroborates our rationale for utilizing window-type data. A larger $l$ has a lower accumulation error. However, it is noteworthy that a model with larger $l$ also costs more training data and time.

\begin{table}[t]
    \centering
    \caption{\textbf{Interpolation error and extrapolation error for sine-Gordon equation.}}
    \begin{tabular}{lcc} 
        \toprule
        \multirow{2}{*}{Models} & Interpolation & Extrapolation \\ 
        & $t \in [1.00, 1.90]$ & $t \in [2.00, 19.90]$ \\
        \midrule
        Vanilla data + recurrent FNO & \textbf{0.0126} & 0.6340 \\
        Vanilla data + full prediction FNO & 0.0274 & blow up \\
        Random data + recurrent FNO & 0.0245 & \textbf{0.0624} \\
        Random data + full prediction FNO & 0.0605 & 0.1928 \\
        \bottomrule
    \end{tabular}
    \label{tab:SG interpolation error and extrapolation error}
\end{table}
\begin{table}[t]
    \centering
    \caption{\textbf{Average relative error for sine-Gordon equation using different window sizes $l$.}}
    \begin{tabular}{lc} 
        \toprule
        Models & Average relative error \\ 
        \midrule
        Random data, recurrent FNO, $l=1$ & 3.1210(blow up)  \\
        Random data, recurrent FNO, $l=2$ & 0.1685  \\
        Random data, recurrent FNO, $l=5$ & 0.0725  \\
        Random data, recurrent FNO, $l=10$ & 0.0606  \\
        Random data, recurrent FNO, $l=20$ & \textbf{0.0326}  \\
        \bottomrule
    \end{tabular}
    \label{tab:SG-l}
\end{table}

\subsection{Klein-Gordon wave equation on irregular domain}\label{subsection: Klein-Gordon wave equation on irregular domain}
In the last PDE, we consider the Klein-Gordon wave equation on an irregular domain. The Klein-Gordon wave equation is closely related to the Schrödinger equation and thereby widely considered in quantum field theory. The examples above are limited to uniform grids on rectangular domains due to the requirement of FFT in FNO. Previous works have successfully extended FNO to general geometries, see our discussion in \autoref{subsection: Neural operator models}. We use the Geo-FNO model proposed in \cite{Li2023Fourier}, which uses a learned deformation mapping the irregular physical space to rectangular computational space. We consider the following ($2+1$)-dimensional Klein-Gordon wave equation with Dirichlet conditions: 
\[
\begin{aligned}
\partial^2_t u - \Delta u + u^3 &= 0, && (t, x) \in (0, 10) \times D, \\
u(t, x) &= u_0(x) , && t \in (0, 10), x \in \partial D, \\
u(t = 0, x) &= u_0(x), && x \in D, \\
\partial_t u(t = 0, x) &= 0 && x \in D. 
\end{aligned}
\]
Here, $D$ is defined by the interior of the polar curves $r = 0.4+0.05(\sin(4\theta) + \cos(3\theta))$, and moving $D$ to the new origin $x = (0.5, 0.5)$. Hence $(0, 1)^2$ is a bounding box of $D$. The data set contains $600$ samples with $100$ testing samples. The spatial-temporal resolution is $2642 \times 200$ with the first snapshot $u_0(x)$ drawn from the Gaussian random field $\mathcal{N}(0, 10^4(-\Delta + 8^2I)^{-6})$ with Neumann condition on $(0, 1)^2$, then restricted back to $D$. $2642$ spatial points include $217$ points on boundary and $2425$ interior points, generated by the algorithm proposed in \cite{Fornberg2015Fast} which is suitable for RBF-FD discretization, see \autoref{fig:KG-plot-domain}. The PDE is solved with the same solver as the sine-Gordon equation. In this example, all Geo-FNO models are FNO-2D. This is because Geo-FNO requires input and query mesh, and we use the scattered nodes on $D$ as mesh. This physical mesh is 2D, so convolution is implemented in a 2D computational latent space. The temporal information is transmitted through the channel dimension. It is also possible to consider a 3D Geo-FNO model duplicating scattered nodes as 3D spatiotemporal points. However, we found that this 3D Geo-FNO has a huge amount of training parameters, takes up a lot of memory, and does not show good numerical performance compared with other models. Hence, we will not include this model here.

We examine the Strichartz estimates on training and validation set, by calculating the ratio of the uniform norm of solution and initial condition on each sample as an estimate for the constant $C$:
\[
C = \frac{\max_{t \in [0, 0.95]} \|u(t, \cdot)\|_\infty}{\|u_0\|_\infty}.
\]
We plot the distribution of $C$ from these $500$ examples in \autoref{fig:KG-plot-estimate} and find that most solutions have a constant sup-norm, which is different from the KdV equation where rough waves appear more frequently. Hence we choose $C = 1$ and clip the output with bound $B = C \cdot \|u_0\|_\infty$ using truncating function $\pi(x)$. \autoref{tab:KG interpolation error and extrapolation error} and \autoref{fig:KG equation} show the performance of different models before and after the clipping operation. We see that clipping does not affect the accuracy and may save the prediction from blowing up in some cases. In \autoref{fig:KG-draw_error_Block-C}, we show the accumulation error of full prediction FNO trained with random data using different values of $C$. We see that $C = 1$ performs the best, which aligns with \autoref{fig:KG-plot-estimate}. \autoref{fig:KG-plot-prediction} shows the clipped prediction by full prediction FNO trained with global random initial data.

\begin{table}[t]
    \centering
    \caption{\textbf{Interpolation error and extrapolation error for Klein-Gordon equation.} Numbers in the brackets indicate the error of the clipped predictions.}
    \begin{tabular}{lccc} 
        \toprule
        \multirow{2}{*}{Models} & Interpolation & \multicolumn{2}{c}{Extrapolation}\\ 
        & $t \in [0.50, 0.95]$ & $t \in [1.00, 2.95]$ & $t \in [3.00, 9.95]$\\
        \midrule
        Vanilla data + recurrent FNO & 0.0197 (0.0197) & 0.1616 (0.1574) & blow up (blow up) \\
        Vanilla data + full prediction FNO & 0.0207 (0.0207) & 0.2832 (0.2831) & 0.3450 (0.3439) \\
        Random data + recurrent FNO & \textbf{0.0167} (0.0168) & 0.1979 (0.1697) & blow up (blow up) \\
        Random data + full prediction FNO & 0.0634 (0.0634) & 0.0896 (\textbf{0.0876}) & blow up (\textbf{0.2327}) \\
        \bottomrule
    \end{tabular}
    \label{tab:KG interpolation error and extrapolation error}
\end{table}

\section{Conclusions and Discussions}\label{section: Conclusions and Discussions}
In this paper, we use the data-driven neural operator method to tackle the long-time prediction problem of nonlinear wave equations. We train the model with temporal length-limited data and predict the complete solution trajectory by iterations. We propose several methods derived from the properties of these equations, which show the potential to reduce the accumulated error and improve algorithm stability:
\begin{itemize}
    \item By the existence and uniqueness of solutions to wave equations, we suggest that the model should take velocity data into account, which can be substituted by the temporal window input.
    \item To approximate the window-to-window mapping, models can learn this mapping directly, making the full prediction in one step, or be embedded into a recurrent network structure and predict step by step. 
    \item Choosing a random initial time of the input can significantly lower the bias of training data and thereby enhance generalization while keeping the size of the training data unchanged. It also helps the model to learn the nonlinear phenomena that appear in long-time evolution.
    \item Conservation quantities of nonlinear wave equations can be used as regularizations adding to the training loss.
    \item Well-posedness based on the Strichartz estimates for the nonlinear wave equation, as a priori upper bound for this set, inspires us to use a clipping operator on the output.
\end{itemize}
We examine these methods by implementing numerical experiments on the KdV equation, sine-Gordon equation, and Klein-Gordon wave equation on an irregular domain. These methods demonstrate effective improvements to the algorithm, yet we acknowledge that the long-time integration problem retains its inherent difficulties. We suggest that the following directions deserve further study in the future: 
\begin{itemize}
    \item Design new network architectures according to inherent features of PDE.
    \item Validate our methods in more scenarios, including different equations, boundary conditions, and geometries. 
    \item Develop algorithms that can accurately model nonlinear phenomena.
    \item Exploit more properties we have not mentioned here from wave equations, such as finite wave speeds, characteristic lines, and analysis related to frequency. 
    \item We use a penalty method as soft constraints on the energy of the output. It is not clear whether a neural operator model can be designed to be energy-preserving, much like many traditional energy-preserving algorithms. 
    \item Combine our approaches with physics-informed models and investigate the numerical effects.
\end{itemize}

\bibliographystyle{plain} 
\bibliography{references.bib}

\appendix
\section{Supplementary Figures}
\begin{figure}[t]
    \centering
    \begin{subfigure}[b]{0.48\textwidth}
        \centering
        \includegraphics[width=\textwidth]{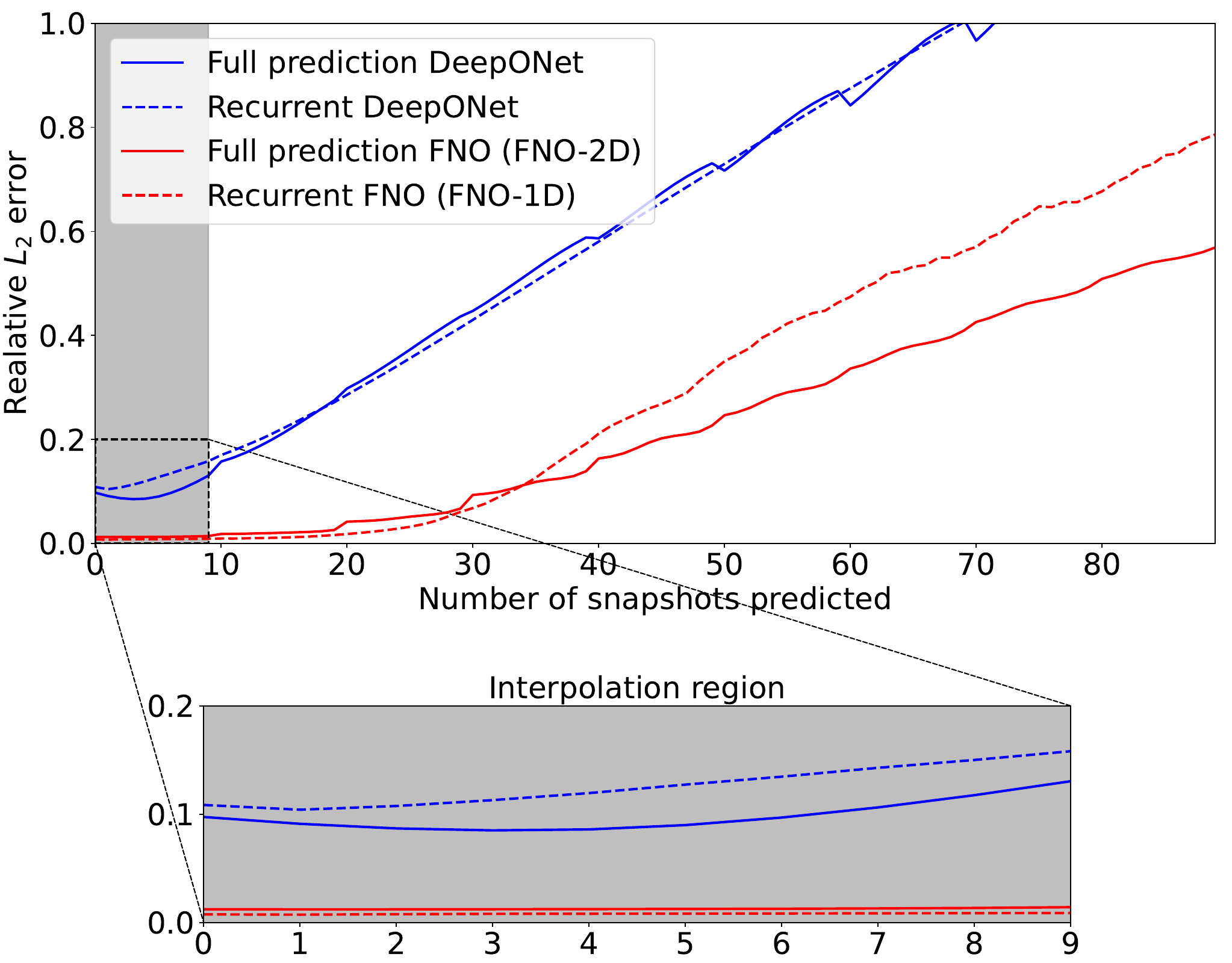}
        \caption{}
        \label{fig:KDV-draw_error_DON-FNO}
    \end{subfigure}

    \vspace{1em}

    \begin{subfigure}[b]{0.48\textwidth}
        \centering
        \includegraphics[width=\textwidth]{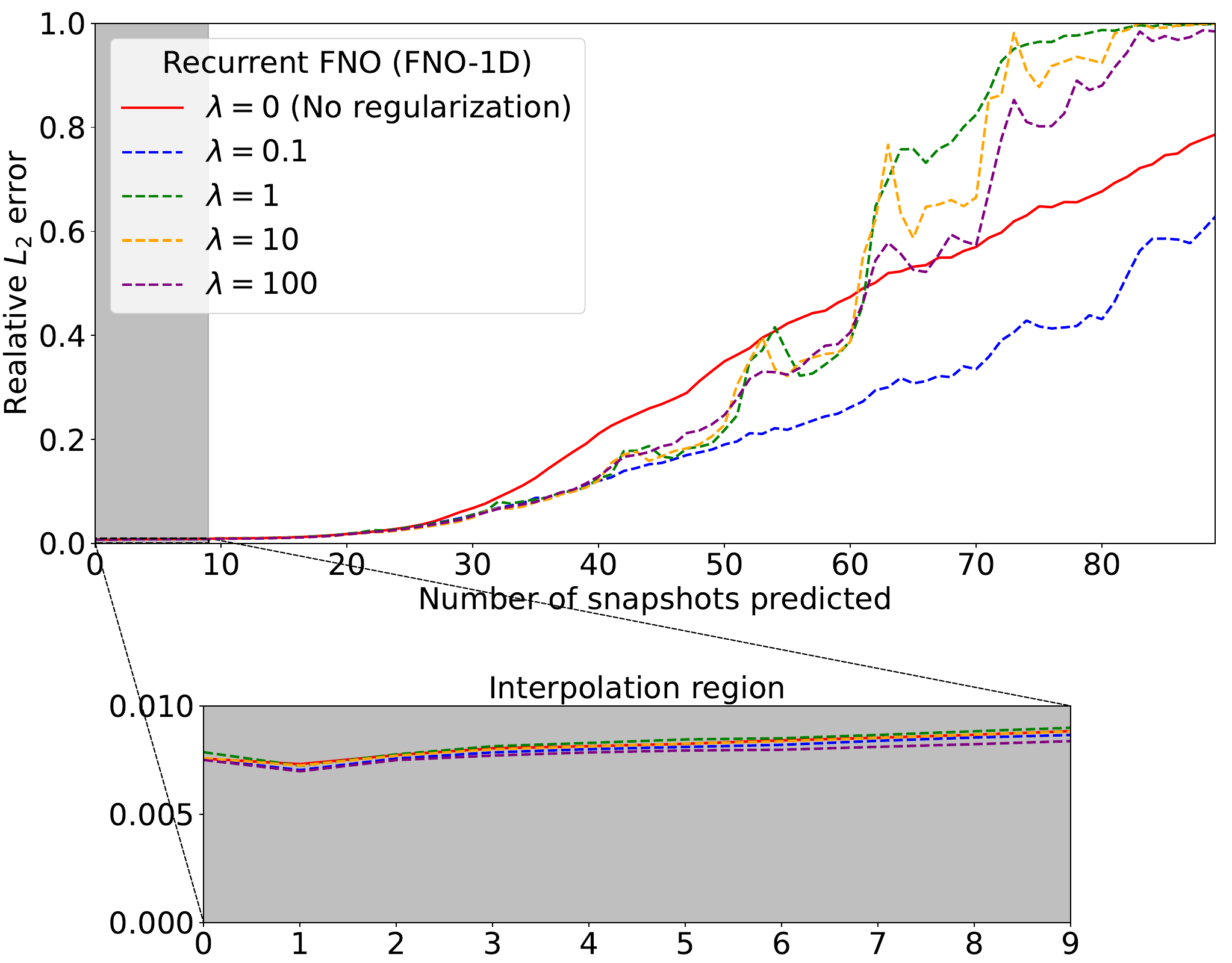}
        \caption{}
        \label{fig:KDV-draw_error_Win-ConLaw}
    \end{subfigure}
    \begin{subfigure}[b]{0.48\textwidth}
        \centering
        \includegraphics[width=\textwidth]{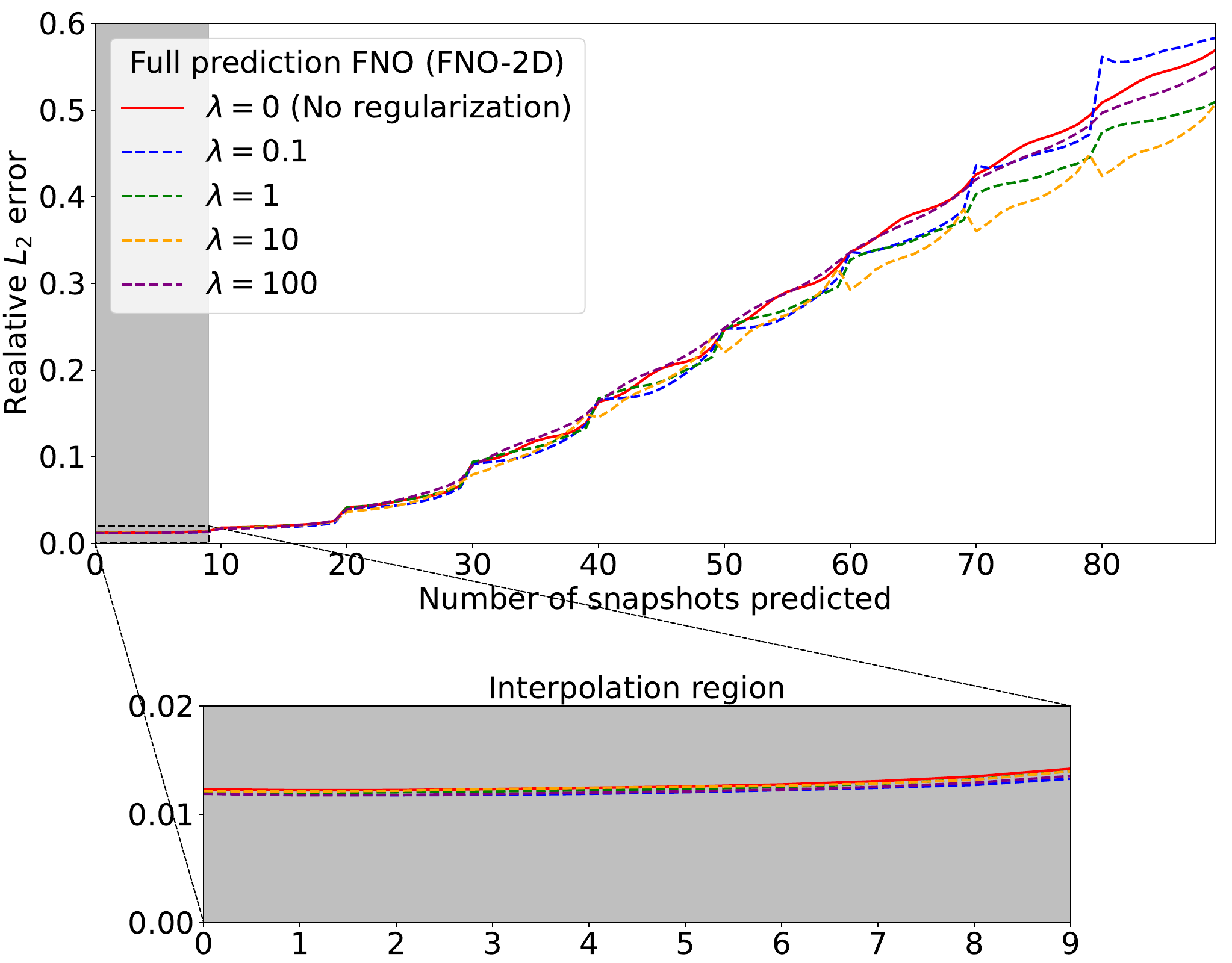}
        \caption{}
        \label{fig:KDV-draw_error_Block-ConLaw}
    \end{subfigure}
    
    \vspace{1em}

    \begin{subfigure}[b]{0.48\textwidth}
        \centering
        \includegraphics[width=\textwidth]{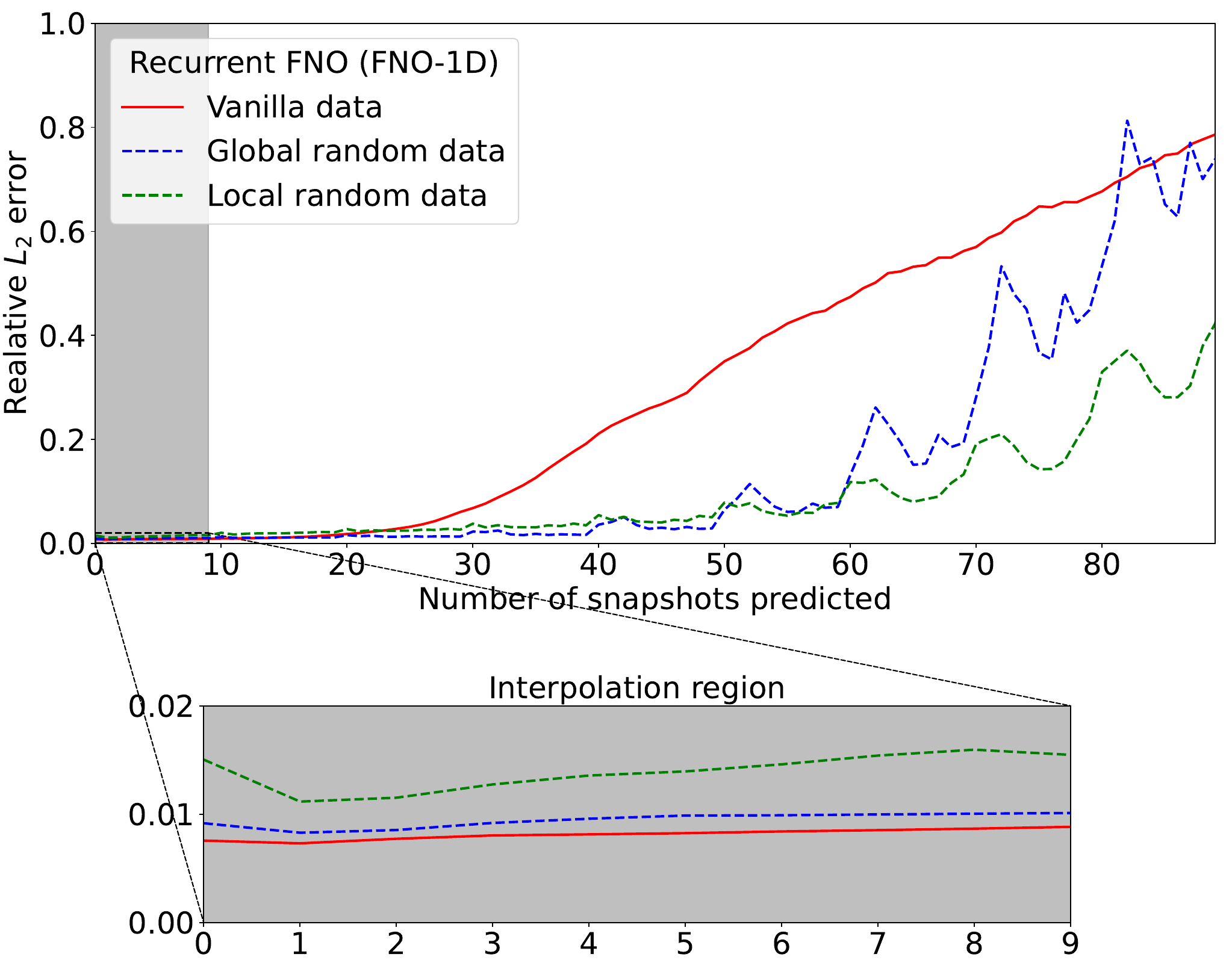}
        \caption{}
        \label{fig:KDV-draw_error_Win-RandomStart}
    \end{subfigure}
    \begin{subfigure}[b]{0.48\textwidth}
        \centering
        \includegraphics[width=\textwidth]{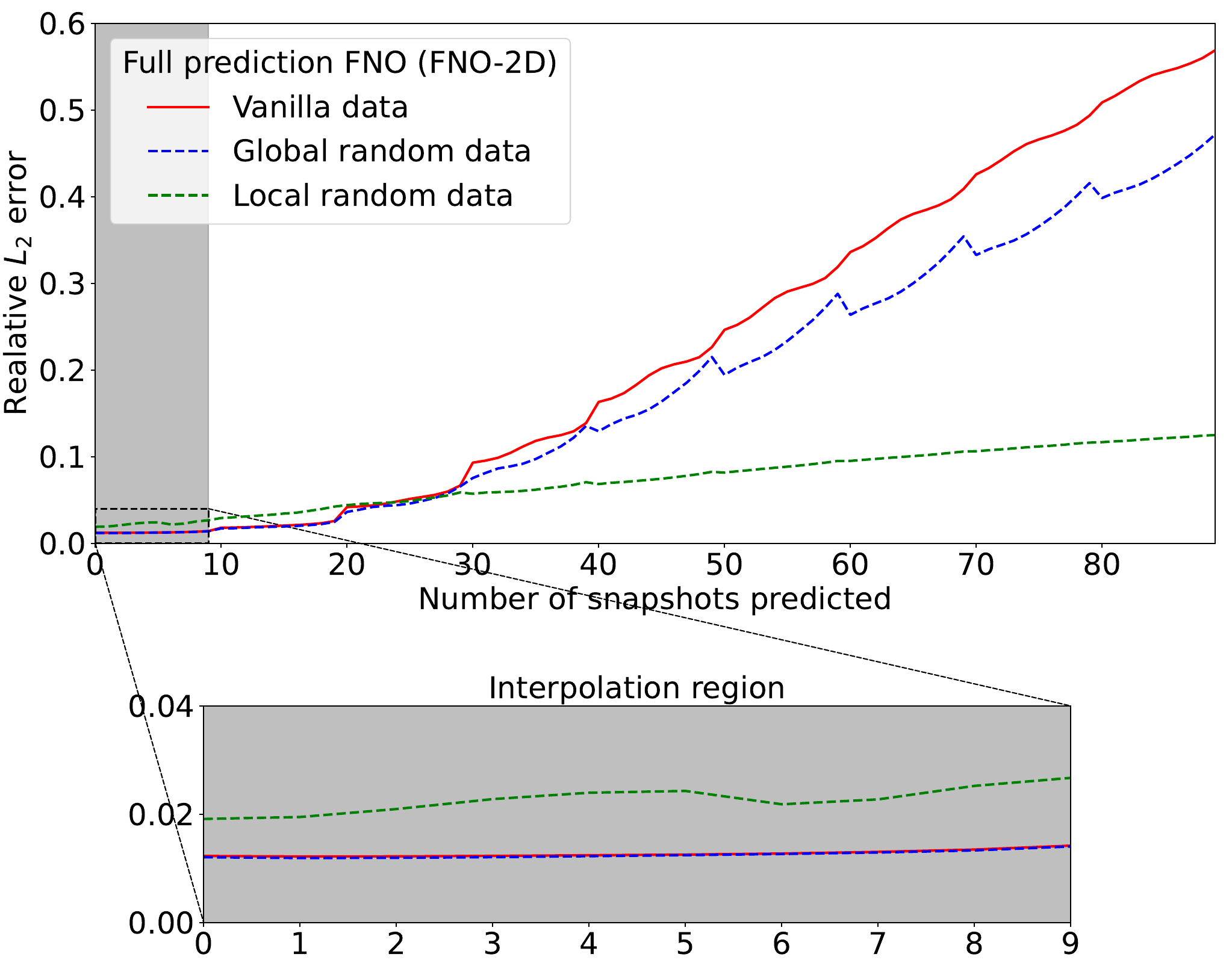}
        \caption{}
        \label{fig:KDV-draw_error_Block-RandomStart}
    \end{subfigure}
    \caption{\textbf{Accumulation error for KdV equation.} (a): Comparison between DeepONet and FNO. (b)(c): Comparison between FNO trained with and without conservation law regularization. (d)(e): Comparison between FNO trained with vanilla data and random data.}
    \label{fig:KdV equation}
\end{figure}

\begin{figure}[t]
    \centering
    \begin{subfigure}[b]{\textwidth}
        \centering
        \includegraphics[width=\textwidth]{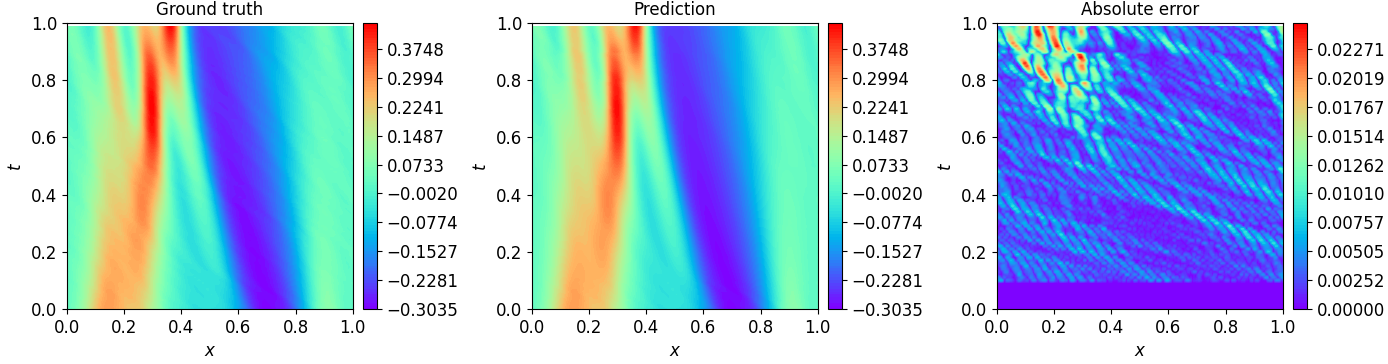}
        \caption{}
        \label{fig:KDV-plot-prediction}
    \end{subfigure}
    
    \vspace{1em}

    \begin{subfigure}[b]{\textwidth}
        \centering
        \includegraphics[width=\textwidth]{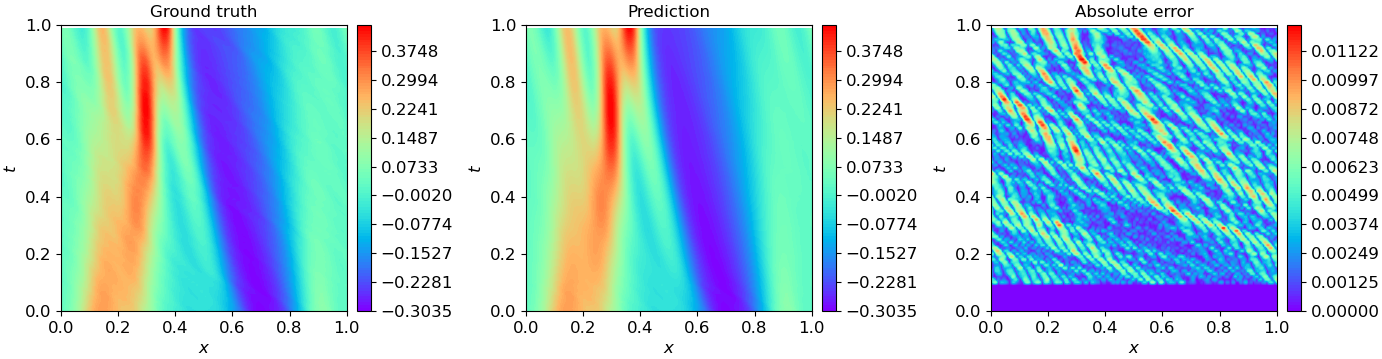}
        \caption{}
        \label{fig:KDV-plot-prediction2}
    \end{subfigure}
    \caption{\textbf{Predictions for KdV equation.} (a): Prediction generated by the full prediction FNO trained with $\lambda = 10$. (b): Prediction generated by the full prediction FNO trained with local random data. The initial $10$ snapshots have zero error as they are the input window.}
\end{figure}

\begin{figure}[t]
    \centering
    \begin{subfigure}[b]{0.48\textwidth}
        \centering
        \includegraphics[width=\textwidth]{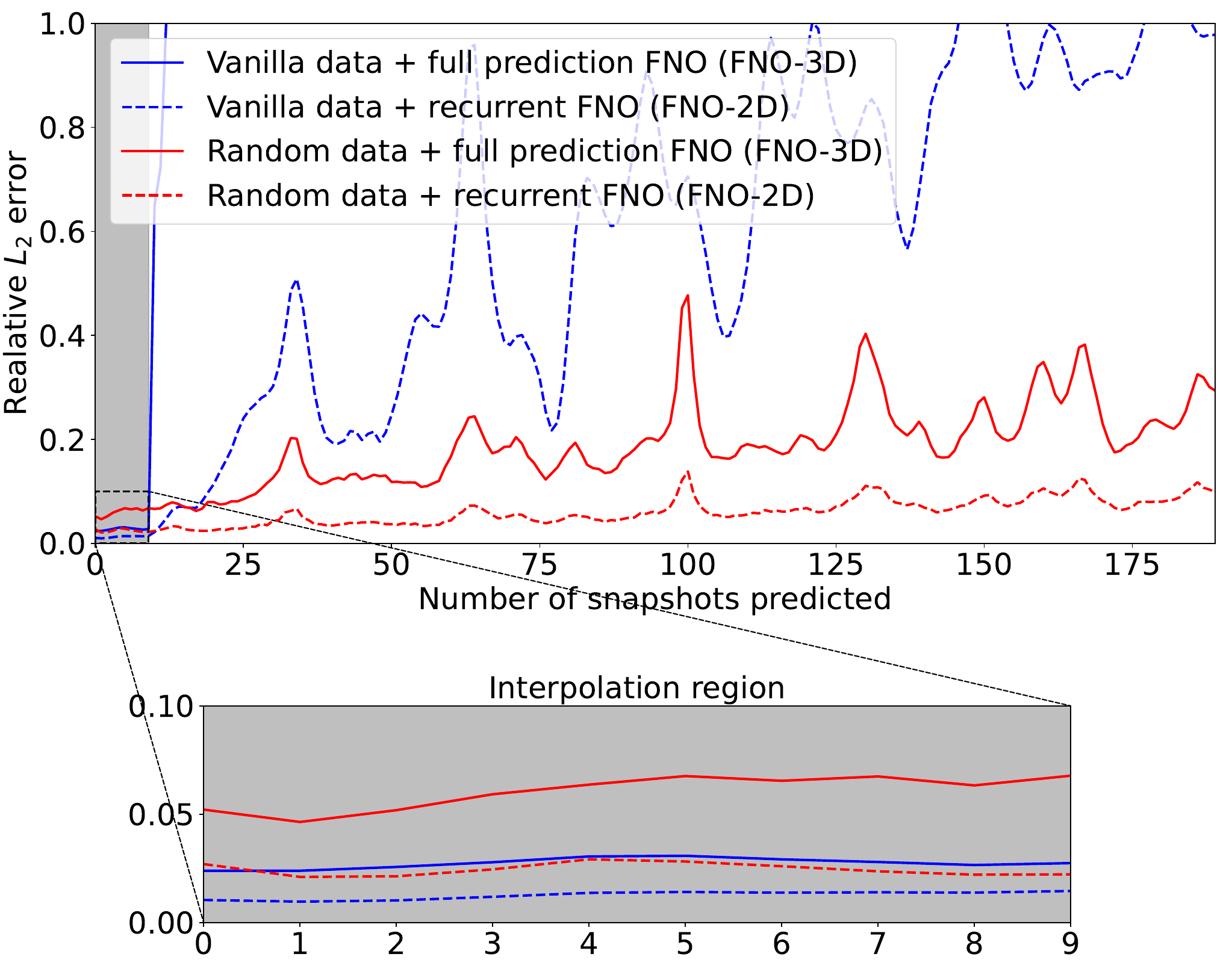}
        \caption{}
        \label{fig:SG-draw_error}
    \end{subfigure}
    \begin{subfigure}[b]{0.48\textwidth}
        \centering
        \includegraphics[width=\textwidth]{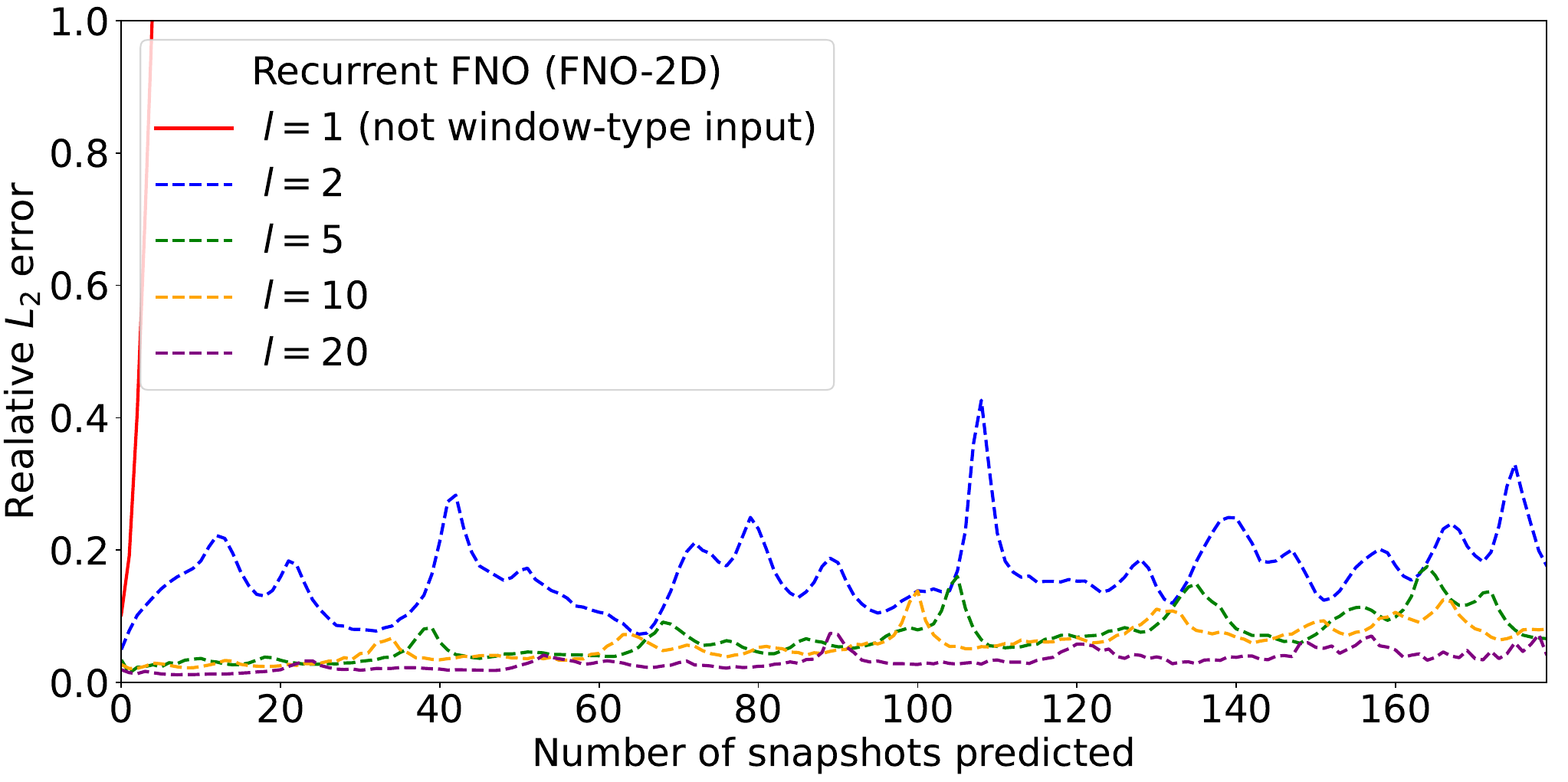}
        \caption{}
        \label{fig:SG-windowsize}
    \end{subfigure}
    \caption{\textbf{Accumulation error for sine-Gordon equation.} (a): Comparison between recurrent and full prediction FNO trained with vanilla data and random data. The grey area shows the trained range and interpolation error of each model. (b): Comparison between recurrent FNO trained with different window sizess $l$.}
\end{figure}

\begin{figure}[t]
    \centering
    \includegraphics[width=0.8\textwidth]{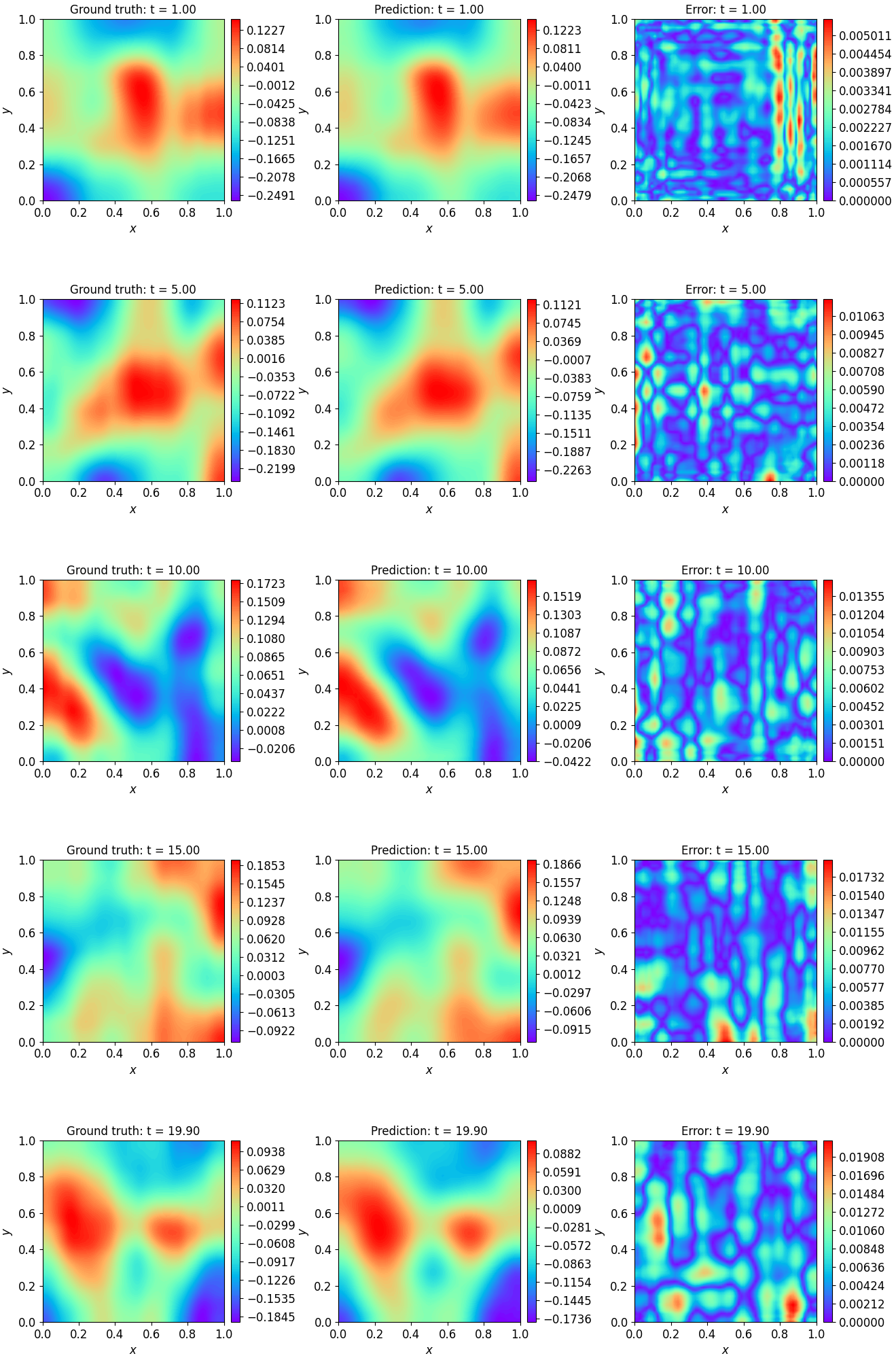}
    \caption{\textbf{Prediction for sine-Gordon equation.} The prediction is made by the recurrent FNO trained with global random initial data.}
    \label{fig:SG-plot-prediction}
\end{figure}

\begin{figure}[t]
    \centering
    \includegraphics[width=0.5\textwidth]{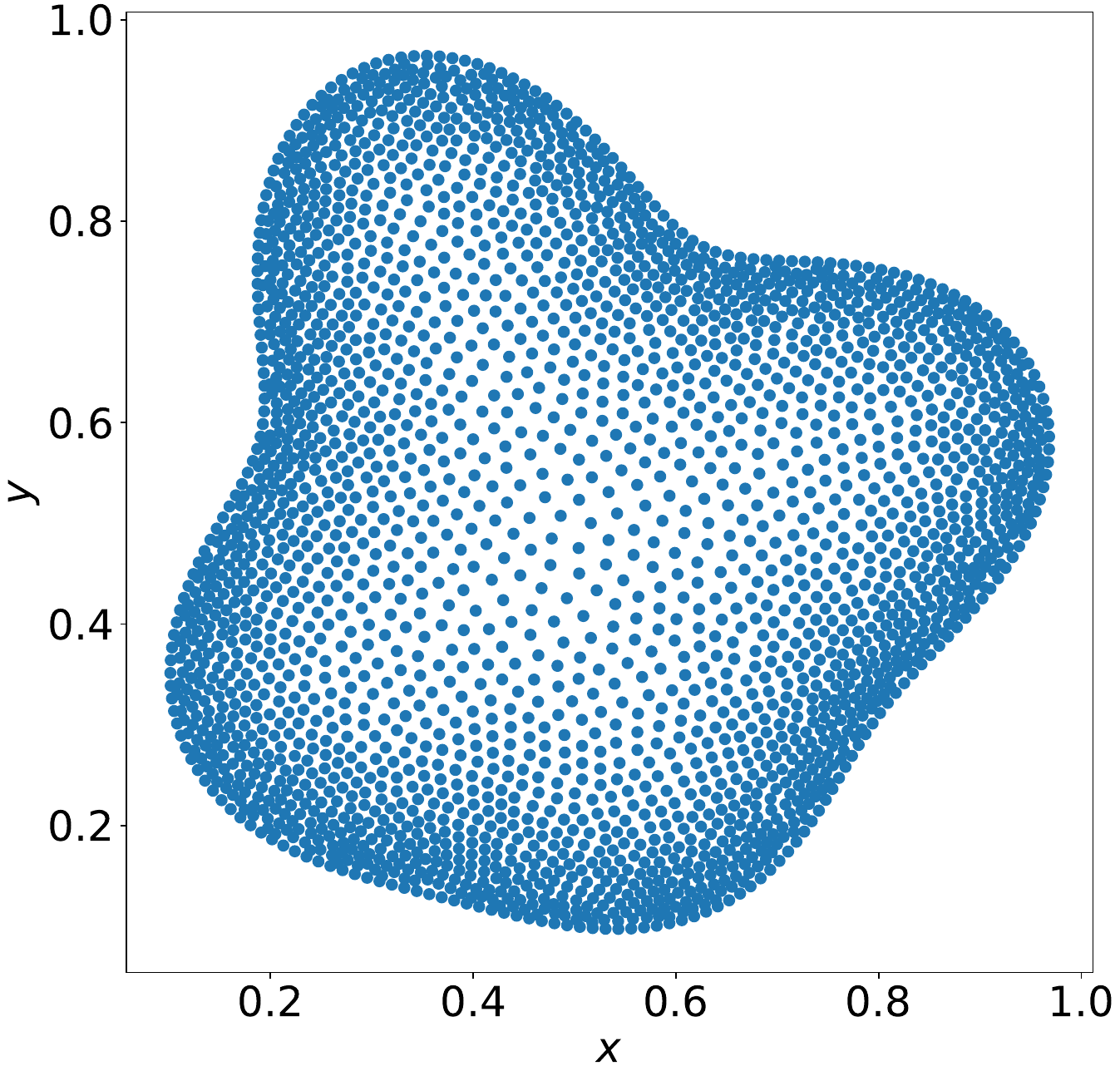}
    \caption{\textbf{Irregular domain $D$ and scattered nodes for Klein-Gordon wave equation.}}
    \label{fig:KG-plot-domain}
\end{figure}
\begin{figure}[t]
    \centering
    \includegraphics[width=0.5\textwidth]{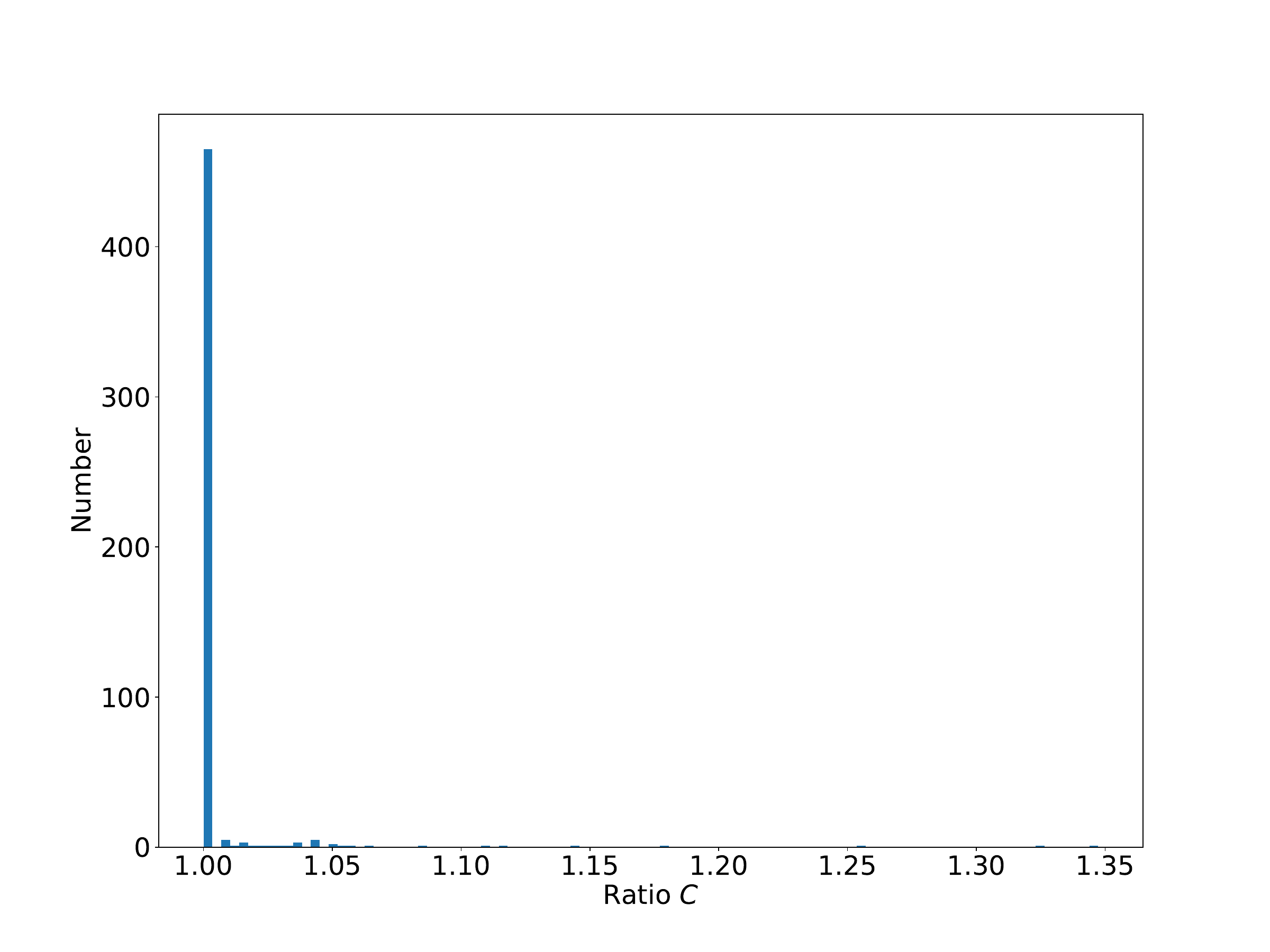}
    \caption{\textbf{Choosing an empirical constant $C$ of Strichartz estimates.}}
    \label{fig:KG-plot-estimate}
\end{figure}
\begin{figure}[t]
    \centering
    \begin{subfigure}[b]{0.48\textwidth}
        \centering
        \includegraphics[width=\textwidth]{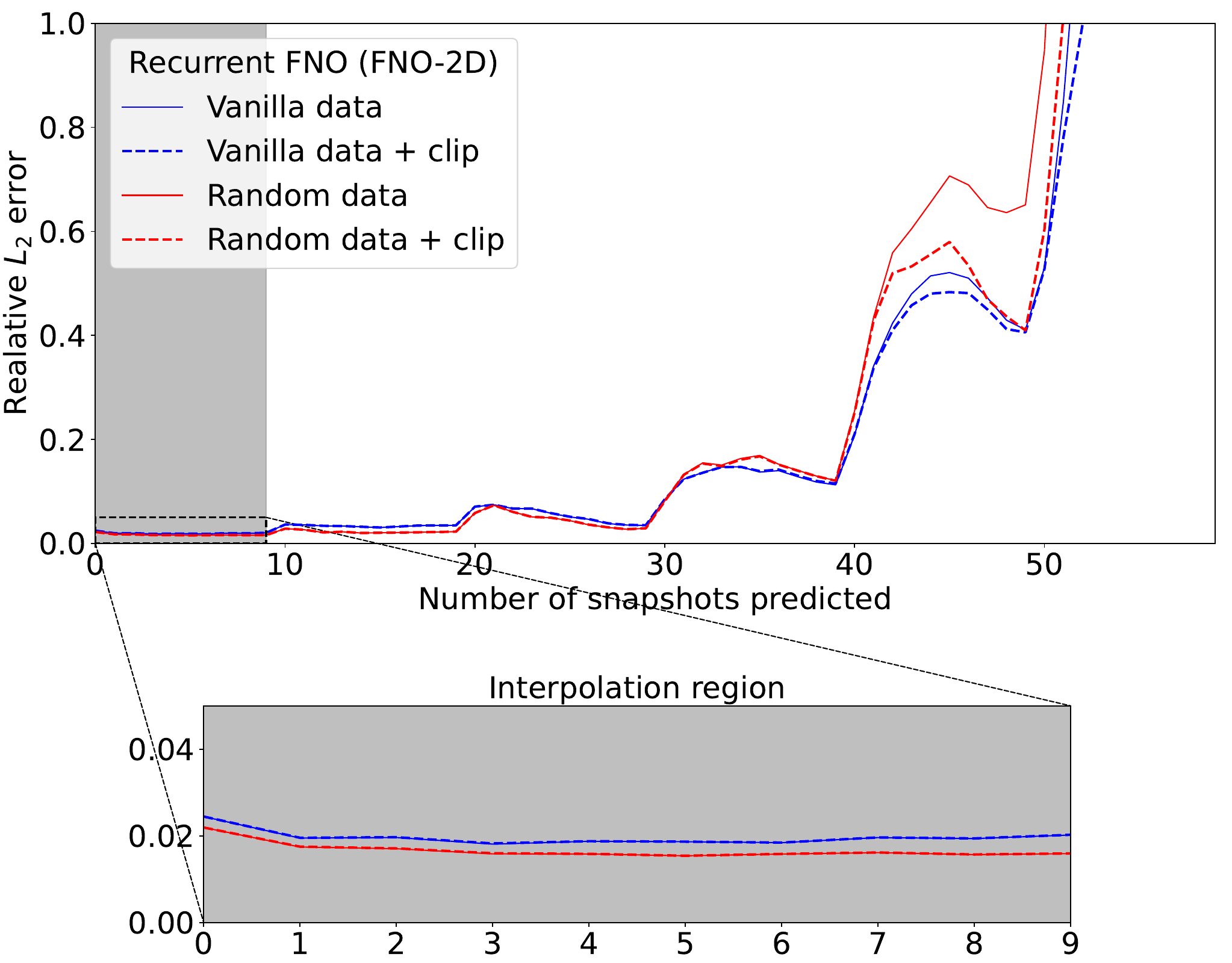}
        \caption{}
        \label{fig:KG-draw_error_Win}
    \end{subfigure}
    \begin{subfigure}[b]{0.48\textwidth}
        \centering
        \includegraphics[width=\textwidth]{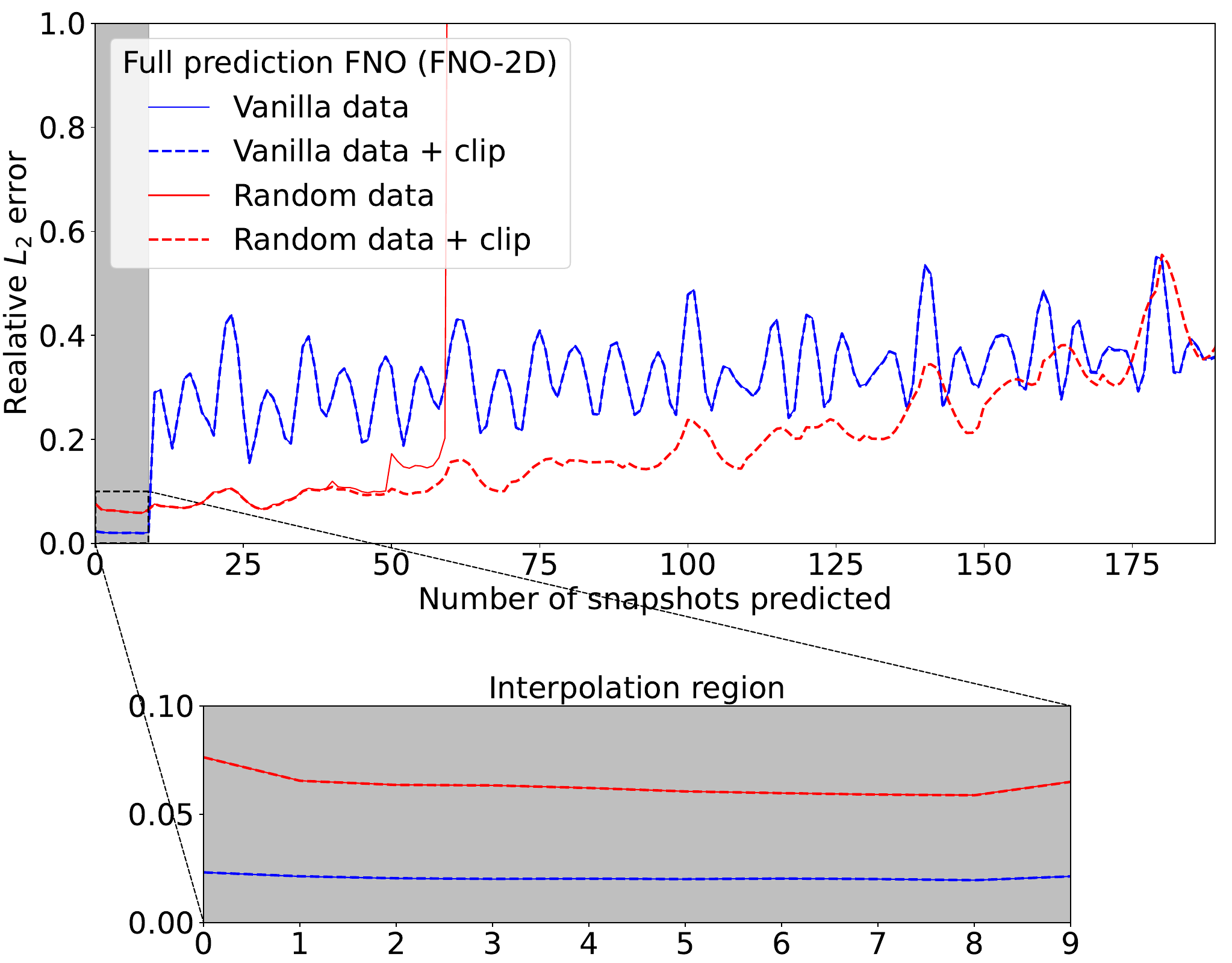}
        \caption{}
        \label{fig:KG-draw_error_Block}
    \end{subfigure}
    \caption{\textbf{Accumulation error for Klein-Gordon equation.} (a): Recurrent FNO. (b): Full prediction FNO. The grey area shows the trained range and interpolation error of each model.}
    \label{fig:KG equation}
\end{figure}
\begin{figure}[t]
    \centering
    \includegraphics[width=0.8\textwidth]{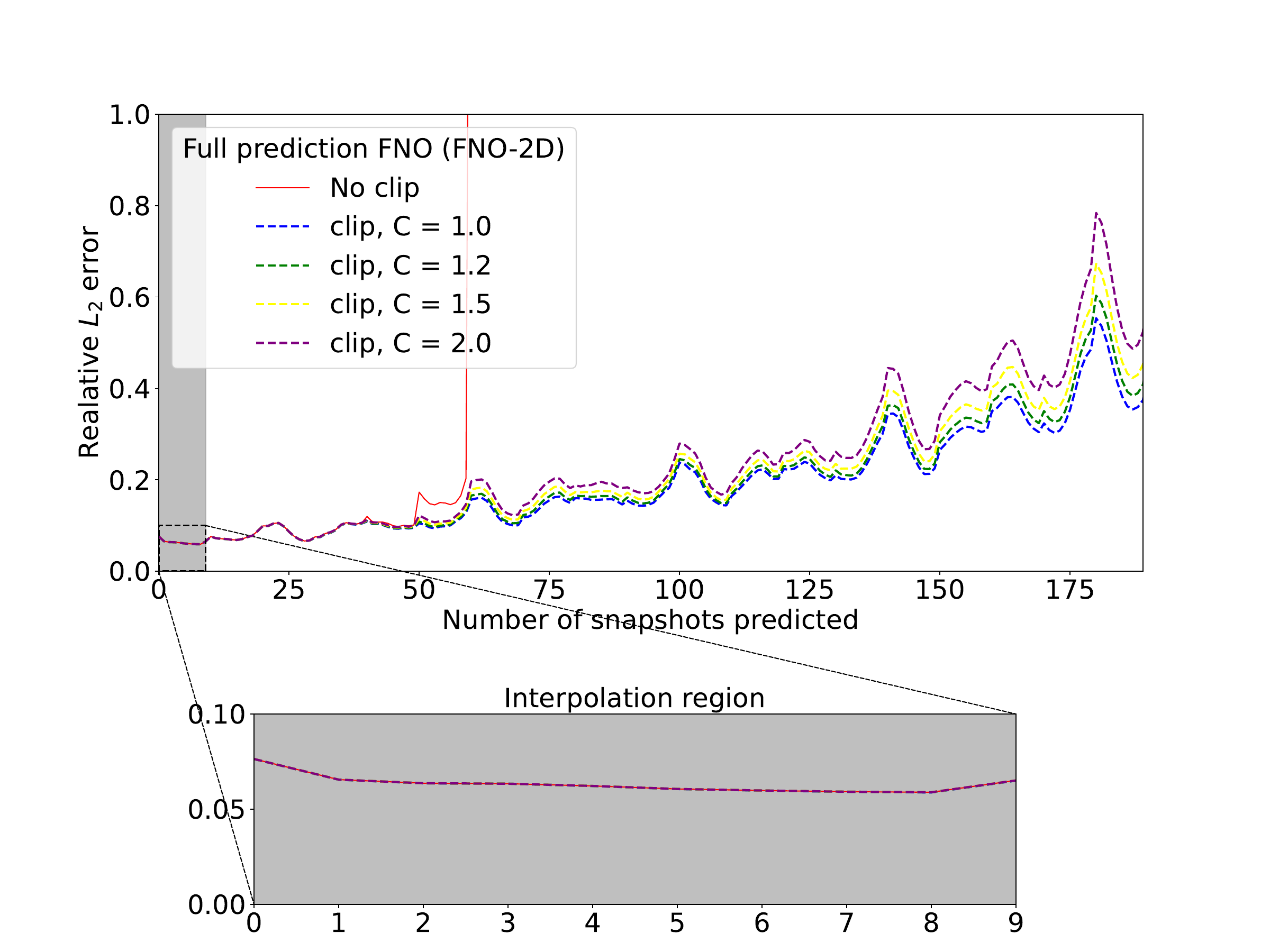}
    \caption{\textbf{Accumulation error with different clip parameters $C$.}}
    \label{fig:KG-draw_error_Block-C}
\end{figure}

\begin{figure}[t]
    \centering
    \includegraphics[width=0.8\textwidth]{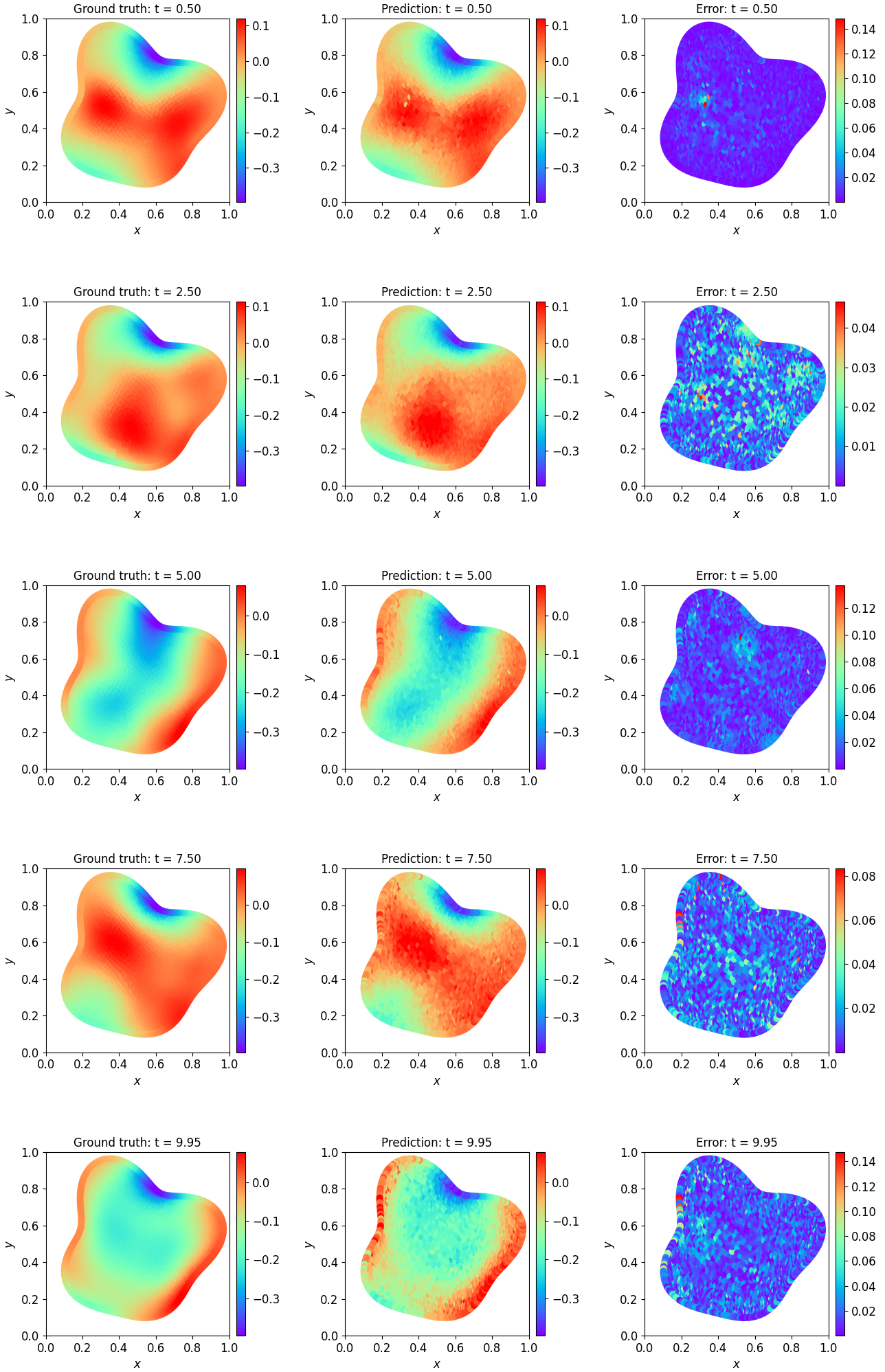}
    \caption{\textbf{Prediction for Klein-Gordon equation on the irregular domain.} The prediction is made by the full prediction FNO trained with global random initial data and clipped to the uniform norm of the initial condition.}
    \label{fig:KG-plot-prediction}
\end{figure}



\end{document}